\documentclass[a4paper,12pt]{amsart}
\usepackage[frenchb,english]{babel}
\usepackage[latin1]{inputenc}
\usepackage{pdfsync}
\usepackage[T1]{fontenc}
\usepackage{amsmath}
\usepackage{amssymb}
\usepackage{amsxtra}
\usepackage{amscd}
\usepackage{amsthm} 
\usepackage{amsfonts}
\usepackage{eucal}
\usepackage[all]{xy}
\usepackage{url}

\newtheorem{thm}{Theorem}[section]
\newtheorem{prop}[thm]{Proposition}
\newtheorem{lem}[thm]{Lemma}
\theoremstyle{definition}
\newtheorem{defi}[thm]{Definition}
\newtheorem{rem}[thm]{Remark}
\numberwithin{equation}{section}
 
\newcommand{\on}{\operatorname}

\newcommand\restr[2]{{ 
  \left.\kern-\nulldelimiterspace    #1  
  \vphantom{\big|}  
  \right|_{#2}  
  }}
\newcommand{\isor}[1]{{\xrightarrow[\raisebox{0.25 em}{\smash{\ensuremath{\sim}}}]{#1}}}

\begin{document}

\title[Shtukas and Langlands correspondence]{
Shtukas for reductive groups and Langlands correspondence for function fields}

\author{Vincent Lafforgue}
\address{Vincent Lafforgue: CNRS et Institut Fourier, UMR 5582, Université Grenoble Alpes, 
 100 rue des Maths, 38610 Gières, France.}
 
\date{\today}
\maketitle

     This text gives an introduction to the Langlands correspondence for function fields and in particular to some recent works in this subject. 
          We begin with a short historical account (all notions used below are recalled in the text).

        The Langlands correspondence  \cite{langlands67} is a conjecture of utmost importance, concerning global fields, i.e. number fields and function fields.    Many excellent surveys are available,  for example  
        \cite{gelbart, bump, jerusalem, taylor-long,   frenkel-lectures, arthur-ICM}. 
        The  Langlands correspondence 
 belongs to  a huge system of  conjectures (Langlands functoriality, Grothendieck's vision of 
 motives, special values of L-functions,  Ramanujan-Petersson conjecture, 
 generalized Riemann hypothesis).  
This system has a remarkable deepness and logical coherence and many cases of these conjectures have already been established.   Moreover the  Langlands correspondence over function fields admits a geometrization, 
  the ``geometric Langlands program'', which is related to conformal field theory in Theoretical Physics.

       Let $G$ be a connected reductive group  over a  global field $F$. 
       For the sake of simplicity we assume $G$ is split. 
                   
    The  Langlands correspondence  relates two fundamental objects, 
    of very different nature, whose definition will be recalled later, 
 \begin{itemize}
 \item the automorphic forms for $G$, 
 \item   the global Langlands parameters , i.e. the conjugacy classes of   morphisms from the  Galois group $\on{Gal}(\overline F/F)$ to the  Langlands dual group $\widehat G(\overline{{\mathbb Q}_\ell})$.   
 \end{itemize}
For $G=GL_{1}$ we have  $\widehat G=GL_{1}$ 
and this is class field theory, which describes the abelianization of  $\on{Gal}(\overline F/F)$ (one particular case of it for  $\mathbb Q$ 
is the  law of quadratic reciprocity, which dates back to Euler, Legendre and  Gauss). 

   Now we restrict ourselves to the case of function fields.      
   
In the case where $G=GL_{r}$ (with $r\geq 2$) the   Langlands correspondence (in both directions)  was proven by  Drinfeld \cite{drinfeld78,Dr1,drinfeld-proof-peterson,drinfeld-compact} for  $r=2$  and by Laurent Lafforgue \cite{laurent-inventiones} for arbitrary  $r$. 
In fact they show the ``automorphic to  Galois'' direction    by using the  cohomology of stacks of shtukas  and the Arthur-Selberg trace formula, 
and deduce from it the reverse direction 
by using the inverse theorems of Weil, Piatetski-Shapiro and Cogdell \cite{inverse-thm} (which are specific to the case of   $GL_{r}$) as well as  Grothendieck's functional equation and Laumon's product formula \cite{laumon-produit} (which are specific to the case of  function fields). 
 Other works using the Arthur-Selberg trace formula for stacks of shtukas are, in chronological order,  Laumon \cite{laumon-rapoport-stuhler,laumon-drinfeld-modular},  Ngo Bao Chau \cite{ngo-jacquet-ye-ulm, ngo-modif-sym,ngo-ngo-elliptique}, Ngo Dac Tuan \cite{ngo-dac-ast, ngo-dac-09},  Lau \cite{eike-lau, eike-lau-duke}, Kazhdan, Varshavsky \cite{var},  Badulescu, Roche \cite{badulescu}.

  In \cite{coh}  we show the ``automorphic to  Galois'' direction   of the Langlands correspondence
         for all   reductive groups over function fields. 
   More precisely we construct a  {\it canonical} decomposition of the 
   vector space of cuspidal automorphic forms, indexed by  global   Langlands parameters. This decomposition is obtained by the  spectral decomposition 
  associated to  the action 
  on this vector space of a  commutative algebra $\mathcal B$ 
  of   ``excursion operators'' such that  each character of   $\mathcal B$ 
  determines a  unique  global Langlands parameter. 
      Unlike previous works, our method is independent on the Arthur-Selberg trace formula. 
     We use the following two ingredients: 
        \begin{itemize}
     \item the classifying stacks of shtukas, introduced by   Drinfeld for $GL_{r}$  
     \cite{drinfeld78,Dr1} and   generalized to all reductive groups and arbitrary number of ``legs'' by  Varshavsky
     \cite{var} (shtukas with several legs were also considered in 
  \cite{ eike-lau, ngo-modif-sym}), 
        \item the geometric  Satake equivalence, due to    Lusztig, Drinfeld, Ginzburg and Mirkovic-Vilonen 
   \cite{hitchin,mv} (it is a fundamental ingredient of the geometric Langlands program, whose idea comes from the fusion of particles in conformal field theory). 
     \end{itemize}
 
In the last sections we  discuss recent works related to the Langlands program over function fields, notably  on the independence on $\ell$ and on the geometric Langlands program. 
  We cannot discuss the works about number fields because there are too many and it is not possible to quote them in this short text.  Let us only mention that, in his lectures at this conference, Peter Scholze will explain local analogues of shtukas over $\mathbb Q_{p}$.

\noindent {\bf Acknowledgements. }
I thank   Jean-Benoît Bost,  Alain Genestier and  Dennis Gaitsgory
 for their crucial help in my research. 
 I am very grateful to  the Centre National de la Recherche Scientifique. 
  The Langlands program  is far from my    first subject 
 and I would never  have been able to devote myself to it without the great freedom 
      given  to CNRS researchers for their works. 
  I thank my colleagues of MAPMO and Institut Fourier for their support. 
I thank  Dennis Gaitsgory for his crucial help in writing the part of this text about geometric Langlands.  I also thank Aurélien Alvarez, Vladimir Drinfeld, Alain Genestier,  Gérard Laumon  and Xinwen Zhu for their help.

 \section{Preliminaries}
 
 \subsection{Basic notions in  algebraic geometry}
 Let $k$ be a field. The ring of functions on the 
 $n$-dimensional affine space $\mathbb A^{n}$ over $k$ is the ring $k[x_{1},...,x_{n}]$ of polynomials in $n$ variables. For any ideal $I$, the quotient $A= k[x_{1},...,x_{n}]/I$
 is the ring of functions on the closed subscheme  of $\mathbb A^{n}$ defined by the equations in $I$ and we obtain in this way all affine schemes (of finite type) over $k$. An affine scheme over $k$ is denoted by $\on{Spec}(A)$ when $A$ is the $k$-algebra of functions on it.   
 It is equipped with the Zariski topology (generated by open subschemes of the form $f\neq 0$ for $f\in A$). 
 It is called a variety when $A$ has no non zero nilpotent element. 
 General schemes and varieties are obtained by gluing. 
 The projective space $\mathbb P^{n}$ over $k$ is the quotient of  $\mathbb A^{n+1}\setminus \{0\}$ by homotheties and can be obtained by gluing $n+1
$ copies of $\mathbb A^{n}$ (which are the quotients of 
$\{(x_{0}, ..., x_{n}), x_{i}\neq 0\}$, for $i=0,...,n$). Closed subschemes (resp. varieties) of  $\mathbb P^{n}$ are called projective schemes (resp varieties) over $k$. Schemes over $k$  have a dimension and a curve is a variety purely  of dimension $1$.

 \subsection{Global fields}   
 A number field is a finite extension of $\mathbb Q$, i.e. a field generated over $\mathbb Q$
by some roots  of a polynomial with  coefficients in $\mathbb Q$.  

A function field $F$  is the field of rational functions 
on an irreducible curve $X$ over a finite field ${{\mathbb F}_q}$. 

We recall that if $q$ is a prime number, ${{\mathbb F}_q}=\mathbb Z/q\mathbb Z$. 
 In general  $q$ is a power of a prime number 
 and all finite fields of cardinal  $q$ are isomorphic to each other
 (although non canonically), hence the   notation   ${{\mathbb F}_q}$.  

The simplest example of a   function field is $F={{\mathbb F}_q}(t)$, 
namely the field of rational functions on the affine line  $X =\mathbb A^{1} $. 
Every   function field is a finite extension of such a field ${{\mathbb F}_q}(t)$. 

Given a function field $F$ there exists a unique smooth projective and geometrically irreducible curve  $X$   over a finite field ${{\mathbb F}_q}$ whose  field of rational functions is $F$: indeed 
for any irreducible curve over ${{\mathbb F}_q}$ we obtain a  smooth projective curve
with the same field of rational functions by resolving the singularities and   adding the points at infinity.  For example   $F={{\mathbb F}_q}(t)$ is the  field of rational functions of the projective line $X=\mathbb P^{1}$  over  ${{\mathbb F}_q}$ (we have added to  $\mathbb A^{1}$ the point  at infinity).   
        
     For the rest of the text we fix    a smooth projective and geometrically irreducible  curve $X$ over 
 ${{\mathbb F}_q}$. We denote by  $F$ the  field of functions of $X$ (but $F$ may also denote a general global field, as in the next subsection).

\subsection{Places of global fields and local fields. }
A place $v$ of a global field $F$ is a non trivial multiplicative norm $F\to \mathbb R_{\geq 0}$, up to equivalence
(where the equivalence relation identifies $\| .\|$ and $\| .\|^{s}$ for any $s>0$). 
The completion $F_{v}$ of the global field $F$  for this norm is called a local field. 
It is a locally compact field and the inclusion $F\subset F_{v}$ determines  $v$. 
Therefore a place is ``a way to complete a global field into a local field''. 

 For any local field there is a canonical normalization of its norm   given
 by the action on its Haar measure. For any non zero element of a global field
 the product of the normalized norms at all places is equal to $1$. 
 
 For example the places of  $\mathbb Q$ are 
\begin{itemize}
\item the archimedean  place, where the completion is  $\mathbb R$ (with normalized norm equal to the usual absolute value),  
 \item for every prime number $p$, the place $p$ where the completion is $\mathbb Q_{p}$ (the normalized norm in  $\mathbb Q_{p}$ of a number $r\in \mathbb Q^{\times}$ 
 is  equal to 
  $p^{-n_{p}(r)}$, where $n_{p}(r)\in \mathbb Z$ is the exponent of  $p$ in the  decomposition of  $r$ as the product of a sign and powers of the prime numbers). 
 \end{itemize}
Thus the local fields obtained by completion of  $\mathbb Q$ are $\mathbb Q_{p}$, for all prime numbers $p$, and $\mathbb R$.
A place $v$ is said to be archimedean if  $F_{v}$ is equal to $\mathbb  R$ or $\mathbb  C$. 
  These places are in finite number for number fields and are  absent for function fields. For each non  archimedean place $v$ we denote by 
 $\mathcal O_{v}$ the ring of integers of  $F_{v}$, consisting of elements of norm $\leq 1$. For example it is $\mathbb Z_{p}$ if $F_{v}=\mathbb Q_{p}$.
 
In the case of  function fields, where we denote by $F$ the field of functions of $X$, the  places are exactly the closed points  of  $X$ (defined as the maximal ideals). 
The closed points  are in bijection with  the orbits under $\on{Gal}(\overline {{\mathbb F}_q}/{{\mathbb F}_q})$ on  $X(\overline{{\mathbb F}_q})$
(Galois groups are recalled below). 
For every closed point $v$ of $X$, we denote by  $n_{v}:F^{\times}\to \mathbb Z$ the  valuation  which associates to a rational function other than $0$  
its vanishing order at $v$. 
We can see  $\mathcal O_{v}$ as the  ${{\mathbb F}_q}$-algebra of functions on the formal neighborhood around  $v$ in $X$ and $F_{v}$ as the  ${{\mathbb F}_q}$-algebra of functions 
on the punctured formal neighborhood. 
We denote  by $\kappa(v)$ the residue field of  $\mathcal O_{v}$; 
it is a finite extension  of ${{\mathbb F}_q}$, whose degree is denoted by $\deg(v)$, therefore it is a finite field with    $q^{\deg(v)}$ elements. The normalized norm on  $F$ associated to  $v$ sends 
$a\in F^{\times}$ to  $q^{-\deg(v)n_{v}(a)}$.

In the example where $X=\mathbb P^{1}=\mathbb A^{1} \cup \infty$,  the unitary  irreducible polynomials in ${{\mathbb F}_q}[t]$ (which is the ring of functions on  $\mathbb A^{1}$) play a role analoguous to that of the prime numbers in $\mathbb Z$: 
 the places of $\mathbb P^{1}$  are 
\begin{itemize}
\item the place $\infty$, at which the completion is  ${{\mathbb F}_q}((t^{-1}))$, 
\item the places associated to unitary  irreducible polynomials in ${{\mathbb F}_q}[t]$  (the degree of such a place is simply the degree of the polynomial). 
  For example the unitary  irreducible polynomial  $t$ corresponds to the point $0\in \mathbb A^{1}$ and the completion at this place is  ${{\mathbb F}_q}((t))$.  
\end{itemize}
We recall that the  local field ${{\mathbb F}_q}((t))$ consists of Laurent series, i.e. sums   $\sum_{n\in \mathbb Z}a_{n}t^{n}$ with  $a_{n}\in {{\mathbb F}_q}$ and $a_{n}=0$ for $n$ negative enough. 

\subsection{Galois groups}
 If $k$ is a field, we denote by  $\overline k$ an algebraic closure of $k$. 
 It is generated over $k$ by the roots of all polynomials  with 
  coefficients in $k$. The separable closure $k^{\on{sep}}\subset \overline k$ consists of the elements whose minimal polynomial over $k$ has a non zero derivative.  We denote by $\on{Gal}(\overline k/k)=\on{Gal}(k^{\on{sep}}/k)$ the group of automorphisms of $\overline k$ (or equivalently of $k^{\on{sep}}$) which act by  the identity on $k$. 
  It is a profinite group, i.e.  a projective limit of finite groups: 
    an  element of  $\on{Gal}(\overline k/k)$ 
    is the same as a family, indexed by the finite Galois extensions $k'\subset \overline k$ of $k$, of elements $\theta_{k'}\in \on{Gal}( k'/k)$,  so that if  $k''\supset k'$, 
       $\restr{\theta_{k''}}{k'}=\theta_{k'}$. 
      We recall that $k'\subset \overline k$ is said to be a finite Galois extension of $k$ if it is a finite dimensional $k$-vector subspace of $k^{\on{sep}}$  and is  stable under the action of 
     $\on{Gal}(\overline k/k)=\on{Gal}(k^{\on{sep}}/k)$ (and then $\on{Gal}(k'/k)$ is a finite group of cardinal equal to the dimension of $k'$ over $k$). 
     
       A simple example is given by finite fields: 
$ \on{Gal}(\overline {{\mathbb F}_q}/{{\mathbb F}_q})$ is equal to the profinite completion  $ \widehat {\mathbb Z}$ of $\mathbb Z$ in such a way that  $1\in \widehat  {\mathbb Z}$ is the Frobenius generator $x\mapsto x^{q}$
(which is an automorphism of $\overline {{\mathbb F}_q}$ equal to identity on  ${{\mathbb F}_q}$). 

We recall that for any ${{\mathbb F}_q}$-algebra, $x\mapsto x^{q}$ is a morphism of ${{\mathbb F}_q}$-algebras, in particular $(x+y)^{q}=x^{q}+y^{q}$. For any scheme $S$ over ${{\mathbb F}_q}$ we denote by  $\on{Frob}_{S}:S\to S$ 
the  morphism acting on functions  by  $\on{Frob}_{S}^{*}(f)=f^{q}$.

We come back to the function field $F$ of  $X$. 
Our main object of interest is the Galois group $ \Gamma=\on{Gal}(\overline F/F)= \on{Gal}(F^{\on{sep}}/F)$.

By the point of view of Grothendieck developed in SGA1, we have an equivalence between
\begin{itemize}
\item the category of finite sets $A$ endowed with a continuous  action of $ \Gamma$
\item the category of  finite separable $F$-algebras
\end{itemize}
where the functor from the first category to the second one maps  $A$ to 
the finite separable $F$-algebra
$ ((F^{\on{sep}})^{A} )^{\Gamma}$  
(here $(F^{\on{sep}})^{A}$ is the direct sum of copies of $F^{\on{sep}}$ indexed by $A$ and $\Gamma$ acts on each copy and permutes them at the same time). 
 We write $  \eta=\on{Spec}(F)$ and $ \overline \eta=\on{Spec}(\overline F)$. 
Then, for any dense open $U\subset X$,  $ \Gamma$ has a profinite quotient 
$\pi_{1}(U, \overline \eta)$ such that a continuous  action of $\Gamma$ on a finite set $A$ factors 
through $\pi_{1}(U, \overline \eta)$ if and only if $\on{Spec}( ((F^{\on{sep}})^{A} )^{\Gamma})$ extends (uniquely) to an étale covering of $U$. 
We will not explain the notion of étale morphism in general and just say that 
a morphism between smooth varieties over a field is  étale if and only if its differential is everywhere invertible. 
 Thus we have an equivalence between 
 \begin{itemize}
\item the category of finite sets $A$ endowed with a continuous  action of $\pi_{1}(U, \overline \eta)$ 
\item the category of  finite étale coverings of $U$. 
\end{itemize}

For any place $v$  the choice of an embedding  $\overline F\subset \overline{F_{v}}$ provides an  inclusion $\on{Gal}(\overline{F_{v}}/F_{v})\subset \on{Gal}(\overline F/F)$  (well defined up to conjugation).  We denote by   $\on{Frob}_{v}\in \on{Gal}(\overline F/F)$  the image 
of any  element of  $ \on{Gal}(\overline F_{v}/F_{v})$
lifting the Frobenius generator  Frobenius $x\mapsto x^{q^{\deg(v)}}$ in    $ \on{Gal}(\overline {\kappa(v)}/\kappa(v))= \widehat {\mathbb Z}$.  
When $U$ is open dense in $X$ as above and $v$ is a place in $U$, the image of 
$\on{Frob}_{v}$ in $\pi_{1}(U, \overline \eta)$ is well defined up to conjugation.

\subsection{A lemma of Drinfeld \cite{drinfeld78}.  }

Let $U\subset X$ open dense as above. For any $i\in I$ we denote by 
$\on{Frob}_{i}$ the ``partial Frobenius'' morphism $U^{I}\to U^{I}$ which sends $(x_{j})_{j\in I}$ to 
$(x'_{j})_{j\in I}$ with $x_{i}'=\on{Frob}_{U}(x_{i})$ and $x_{j}'=x_{j}$ for $j\neq i$. 
For any scheme $T$ and any morphism $T\to U^{I}$, we say that a morphism 
$a:T\to T$ is ``above'' $\on{Frob}_{i}$ if the square 
$$\xymatrix{ 
T \ar[d] \ar[r]^{a} & T  \ar[d]  \\
U^{I} \ar[r]^{\on{Frob}_{i}} & U^{I} }$$
is commutative.

 \begin{lem}     \label{lem-Frob-partiels-drinfeld-ensembles}
  We have an equivalence of categories between 
    \begin{itemize}
     \item the category of finite sets $A$ endowed with a continuous  action of 
  $  (\pi_{1}(U, \overline \eta))^{I}$, 
  \item  the category  of finite étale coverings $T$ of  $U^{I}$, 
  equipped with partial Frobenius morphisms, i.e. morphisms $F_{\{i\}}$ above $\on{Frob}_{i}$, commuting with each other, and whose composition is  $\on{Frob}_{T}$.  
    \end{itemize}
      \end{lem}
      
      The functor from the first category to the second one is the following: 
      if the action of $(\pi_{1}(U, \overline \eta))^{I}$ on $A$ factorizes through
      $\prod _{i\in I}\on{Gal}(U_{i}/U)$ where for each $i$, 
      $U_{i}$ is a finite étale Galois covering of $U$ (and $\on{Gal}(U_{i}/U)$ is its automorphism group), then the image by the functor  is  $(\prod_{i\in I}U_{i})\times_{\prod _{i\in I}\on{Gal}(U_{i}/U)} A$,  
      equipped with the partial  Frobenius morphisms 
      $F_{\{i\}}$ given by 
      $\big(\on{Frob}_{U_{i}}\times \prod_{j\neq i }\on{Id}_{U_{j}}\big)\times \on{Id}_{A}$.

 \subsection{Split connected reductive groups and bundles}

We denote by $\mathbb G_{m}=GL_{1}$ the multiplicative group. 
 A split torus over a field $k$ is an algebraic group $T$ which  is isomorphic to $\mathbb G_{m}^{r}$ for some $r$. 
 
 A connected reductive group over a field $k$ is a connected, smooth,  affine algebraic group whose extension to $\overline k$ has a trivial unipotent radical (i.e. any  normal, smooth,  connected,  unipotent  subgroup scheme of it is trivial). A connected reductive group $G$ over   $k$ is said to be split if it has a split maximal torus $T$. Then (after chosing a Borel subgroup containing $T$) the lattices $\on{Hom}(\mathbb G_{m}, T)$ and $\on{Hom}(T,\mathbb G_{m})$ are called the coweight and weight lattices of $G$. 
 The split connected reductive groups over a field $k$ are exactly the quotients by central finite subgroup schemes of products of $\mathbb G_{m}$, simply-connected split  groups in the four series $SL_{n+1}$, $Spin_{2n+1}$, $Sp_{2n}$, $Spin_{2n}$, and five simply-connected split exceptional groups.

 Let $G$ be a split connected reductive group over a field $k$, and 
 $X$  a scheme over $k$. Then a $G$-bundle over $X$ is a morphism $Y\to X$, together which an action of $G$ on the fibers which is simply transitive. A $GL_{r}$-bundle $\mathcal E$ gives rise to the vector bundle of rank $r$ equal to $\mathcal E\times_{GL_{r}}\mathbb A^{r}$ and the notions are equivalent.

\section{Reminder on automorphic forms}

For the moment we take  $G=GL_{r}$. When the global field is  $\mathbb Q$, an automorphic form  (without level at finite places) is a  function on the  quotient $GL_{r}(\mathbb Z)\backslash GL_{r}(\mathbb R)$ (the best known example is the particular case of modular functions, for which  $r=2$). 
This quotient classifies the free  
 $\mathbb Z$-modules  (or, equivalently, projective $\mathbb Z$-modules)  $M$ of rank $r$  equipped with a trivialization $M\otimes_{\mathbb Z} \mathbb R=\mathbb R^{r}$ (i.e.  an embedding of  $M$ as a lattice in $\mathbb R^{r}$).  Indeed if we choose a basis of  $M$ over $\mathbb Z$ its embedding in  $\mathbb R^{r}$ is given by a matrix in  $GL_{r}(\mathbb R)$ and the change of the basis of $M$  gives the  quotient by $GL_{r}(\mathbb Z)$. 
 
 Now we come back to our function field $F$. To explain the analogy with  $\mathbb Q$ we choose a  place $v$ of  $X$ (of degree $1$ to simplify)   playing the role of the archimedean place of  $\mathbb Q$ (but this choice is not natural and will be forgotten in five lines). 
 An analogue of a projective $\mathbb Z$-module   $M$ of rank $r$  equipped with a trivialization $M\otimes_{\mathbb Z} \mathbb R=\mathbb R^{r}$ is a  vector bundle of rank 
  $r$ over  
$X$ equipped with a trivialization on the formal neighborhood around $v$. 
 Now we forget the trivialization on the formal neighborhood around $v$ 
 (because we do not want  to introduce a level at $v$) 
 and then we forget the choice of  $v$.  
 
Thus an automorphic form (without level at any place) for  $GL_{r}$ 
is  a function on  the set 
 $\on{Bun}_{GL_{r}}({{\mathbb F}_q})$ of  isomorphism classes of vector bundles 
 of rank   $r$ over $X$.

  Now we consider the case of a general group $G$.       
   From now on we denote by   $G$ a connected reductive group over $F$,    assumed to be split for simplicity. 
      An automorphic form (without level) for  $G$ is a  function on  the set $\on{Bun}_{G}({{\mathbb F}_q})$ of isomorphism classes of $G$-bundles
     over $X$.

   \begin{rem}\label{rem-artin} This remark can be skipped.   
In fact the  $G$-bundles over $X$  have finite automorphism groups.  
Therefore it is more natural to consider  $\on{Bun}_{G}({{\mathbb F}_q})$ as a  groupoid, i.e.  a category where all arrows are  invertible. 
 It is the groupoid of points over  ${{\mathbb F}_q}$ of the Artin  stack $\on{Bun}_{G}$ over  ${{\mathbb F}_q}$ whose groupoid of   $S$-points
(with $S$ a scheme over ${{\mathbb F}_q}$) classifies the $G$-bundles    over  
$X\times S$.  We refer to  \cite{laumon-moret-bailly} for the notion of Artin stack  
and we say only that  examples of  Artin stacks are  given by  the quotients of  algebraic varieties by   algebraic groups. In Artin stacks the   automorphism groups of the points are 
algebraic groups (for example in the case of a quotient they are the stabilizers of  the action).  
  The Quot construction  of  Grothendieck implies that 
  $\on{Bun}_{G}$ is   an Artin stack  
      (locally it is even the quotient of a smooth algebraic variety by a smooth algebraic group).  The automorphisms groups of points in the groupoid 
      $\on{Bun}_{G}({{\mathbb F}_q})$ are finite, in fact they are the   points over   ${{\mathbb F}_q}$ 
      of automorphisms groups of points in $\on{Bun}_{G}$, which are 
      algebraic groups of finite type. 
     \end{rem}

  It is convenient to impose a  condition relative to the center $Z$ of $G$.  
    From now on   we fix a subgroup  $\Xi$ of finite index   in  $\on{Bun}_{Z}({{\mathbb F}_q})$ (for example the trivial subgroup if  $Z$ is finite) 
      and we consider functions on    $\on{Bun}_{G}({{\mathbb F}_q})/\Xi$.    
However, except when $G$ is a torus,  $\on{Bun}_{G}({{\mathbb F}_q})/\Xi$ is still infinite. 
 To obtain  vector spaces of  finite dimension we now restrict ourselves to cuspidal automorphic forms.  
 
      For any field $E\supset \mathbb Q$, we denote by  $C_{c}^{\mathrm{cusp}}(\on{Bun}_{G}({{\mathbb F}_q})/\Xi,E) $ the $E$-vector space of finite dimension consisting of    {\it cuspidal}   functions   on $\on{Bun}_{G}({{\mathbb F}_q})/\Xi$. It is defined as the intersection of the kernel of all ``constant term'' morphisms $C_{c} (\on{Bun}_{G}({{\mathbb F}_q})/\Xi,E)\to C(\on{Bun}_{M}({{\mathbb F}_q})/\Xi,E) $
      (which are given by the correspondence  $\on{Bun}_{G}({{\mathbb F}_q})\leftarrow \on{Bun}_{P}({{\mathbb F}_q})\rightarrow \on{Bun}_{M}({{\mathbb F}_q})$ and involve only finite sums), 
     for all proper parabolic subgroups   $P$ of  $G$ with associated Levi quotient  $M$ (defined as the quotient of $P$ by its unipotent radical).  For readers who do not know these notions, we recall that in the case of  $GL_{r}$ a parabolic subgroup  $P$ is conjugated to a subgroup of upper  block triangular matrices and that the associated    Levi quotient $M$ is isomorphic  
   to the group of   block diagonal matrices.   
 It is legitimate in the Langlands correspondence to restrict oneself  to 
     cuspidal automorphic forms  because all automorphic forms for  $G$ 
     can be understood from  
     cuspidal automorphic forms  for  $G$ and for the Levi quotients of its parabolic subgroups.

     Let   $\ell$ be a prime number not dividing $q$. 
    To simplify the notations we assume that  ${\mathbb Q}_\ell$ contains a square root of  $q$
    (otherwise replace ${\mathbb Q}_\ell$ everywhere by a finite extension containing a square root of  $q$). 
 For Galois representations we have to   work with  coefficients in ${\mathbb Q}_\ell$  and $\overline{{\mathbb Q}_\ell}$, and not  $\mathbb Q, \overline{\mathbb Q}$ and even ${\mathbb C}$ (to which  $\overline{{\mathbb Q}_\ell}$  is isomorphic algebraically but not topologically) 
      because the Galois representations which are continuous with coefficients in ${\mathbb C}$ always have a finite image (unlike  those with coefficients in  $\overline{{\mathbb Q}_\ell}$) 
      and  are  not enough to match automorphic forms in  the  Langlands correspondence.  Therefore, even if the notion of cuspidal  automorphic form
      is (in our case of function fields) algebraic,  to study the Langlands correspondence we will  consider cuspidal  automorphic forms with coefficients in $E={\mathbb Q}_\ell$ or $\overline{{\mathbb Q}_\ell}$.

\section{Class field theory for function fields} \label{CFT}

It was developped by  Rosenlicht and Lang \cite{livre-serre}. Here we consider only the unramified case. 

Let $\on{Pic}$ be the relative Picard scheme of $X$ over ${{\mathbb F}_q}$, 
whose definition is that,  for any scheme $S$ over ${{\mathbb F}_q}$, $\on{Pic}(S)$ 
(the set of morphisms $S\to \on{Pic}$) 
 classifies the isomorphism classes $[\mathcal E]$ of  line bundles $\mathcal E$ on $X\times S$ (a line bundle is a vector bundle of rank $1$, so it is the same as a $\mathbb G_{m}$-bundle). The relation with $\on{Bun}_{GL_{1}}$ is that 
 $\on{Bun}_{GL_{1}}$ can be identified with the quotient of $\on{Pic}$ by the trivial action of $\mathbb G_{m}$.

Let $\on{Pic}^{0}$ be the neutral component of  $\on{Pic}$, 
i.e. the kernel of the degree morphism $\on{Pic}\to \mathbb Z$. It is an abelian variety over ${{\mathbb F}_q}$, also called the   jacobian of $X$. 

Class field theory states 
(in the unramified case to which we restrict ourselves in this text) that 
there is a canonical isomorphism  
\begin{gather}\label{isom-CFT}\big(\pi_{1}(X,\overline \eta)\big)^{\on{ab}} \times_{ \widehat {\mathbb Z} }\mathbb Z\overset {\thicksim}{\to} \on{Pic}({{\mathbb F}_q}) \end{gather}
 characterized by the fact that for any place $v$ of $X$, it sends  $\on{Frob}_{v}$  to $[\mathcal O(v)]$, where $\mathcal O(v)$ is the line bundle on $X$ whose sections are the functions on $X$ with a possible pole of order $\leq 1$ at $v$. 

The isomorphism \eqref{isom-CFT} implies that for any $a\in \on{Pic}({{\mathbb F}_q})$ of non zero degree we can associate to  any (multiplicative) character $\chi$ of the finite abelian group $\on{Pic}({{\mathbb F}_q})/a^{\mathbb Z}$ (with values in any field, e.g. 
$\overline{{\mathbb Q}_\ell}$ for $\ell$  prime to $q$) a character 
$\sigma(\chi)$  of $\pi_{1}(X,\overline \eta)$.
We now give a geometric construction of $\sigma(\chi)$, which is in fact  the key step in the proof of the isomorphism \eqref{isom-CFT}. 
 
The Lang isogeny $L:\on{Pic}\to \on{Pic}^{0}$ is such that, for any scheme $S$ over ${{\mathbb F}_q}$  and every line bundle $\mathcal E$ on $X\times S$, 
$[\mathcal E]\in \on{Pic}(S)$ is sent by $L$ to 
$[\mathcal E^{-1}\otimes (\on{Frob}_{S}\times \on{Id}_{X})^{*}(\mathcal E)]\in\on{Pic}^{0}(S)$. We note that  $[(\on{Frob}_{S}\times \on{Id}_{X})^{*}(\mathcal E)]\in \on{Pic}(S)$   is  the image by $\on{Frob}_{\on{Pic}}$ of $[\mathcal E]\in \on{Pic}(S)$. 
The Lang isogeny is surjective and its kernel is $\on{Pic}({{\mathbb F}_q})$. 
For any finite set $I$ and any family $(n_{i})_{i\in I}\in \mathbb Z^{I}$ satisfying $\sum_{i\in I}n_{i}=0$, we consider the Abel-Jacobi morphism $AJ:X^{I}\to \on{Pic}^{0}$
sending $(x_{i})_{i\in I}$ to the line bundle $\mathcal O(\sum_{i\in I }n_{i}x_{i})$. 
We form the fiber  product
$$\xymatrix{\ar @{} [dr] |{\square}
\on{Cht}_{I,(n_{i})_{i\in I}} \ar[d]^{\mathfrak p} \ar[r] & \on{Pic}  \ar[d]^{L} \\
X^{I} \ar[r]_{AJ} & \on{Pic}^{0} }$$
and see that $\mathfrak p$ is a Galois covering of $X^{I}$ with Galois group 
$\on{Pic}({{\mathbb F}_q})$. 
Thus, up to an automorphism group ${{\mathbb F}_q}^{\times}$  which we neglect, for any scheme $S$ over ${{\mathbb F}_q}$, 
$\on{Cht}_{I,(n_{i})_{i\in I}}(S)$ classifies 
\begin{itemize}
\item 
morphisms $x_{i}:S\to X$ 
\item a line bundle $\mathcal E$ on $X\times S$ 
 \item  an isomorphism  $\mathcal E^{-1}\otimes (\on{Frob}_{S}\times \on{Id}_{X})^{*}(\mathcal E)\simeq \mathcal O(\sum_{i\in I }n_{i}x_{i})$.
 \end{itemize}
  Moreover $\on{Cht}_{I,(n_{i})_{i\in I}}  $ is equipped with partial Frobenius morphisms $F_{\{i\}}$ sending $\mathcal E$ to  $\mathcal E\otimes \mathcal O(n_{i}x_{i})$. The morphism 
$F_{\{i\}}$ is above $\on{Frob}_{i}:X^{I}\to X^{I}$, because $(\on{Frob}_{S}\times \on{Id}_{X})^{*}(\mathcal O( x_{i}))=\mathcal O( \on{Frob}_{S}(x_{i}))$. 
Taking the quotient by $a^{\mathbb Z}$ we obtain a {\it finite} Galois covering
$$\xymatrix{ 
\on{Cht}_{I,(n_{i})_{i\in I}}/a^{\mathbb Z} \ar[d]^{\mathfrak p}   \\
X^{I}  }$$
with Galois group $\on{Pic}({{\mathbb F}_q})/a^{\mathbb Z}$ and equipped with the 
  partial Frobenius morphisms $F_{\{i\}}$.  
Then Drinfeld's lemma gives rise to a morphism 
$\alpha_{I,(n_{i})_{i\in I}}:\pi_{1}(X,\overline \eta)^{I}\to \on{Pic}({{\mathbb F}_q})/a^{\mathbb Z}$. 
The character $\sigma(\chi)$ of $\pi_{1}(X,\overline \eta)$ is characterized by the fact that for any $I$ and $(n_{i})_{i\in I}$  with sum $0$,  $\chi\circ \alpha_{I,(n_{i})_{i\in I}}=\boxtimes _{i\in I} \sigma(\chi)^{n_{i}}$ and this gives in fact a construction of $\sigma(\chi)$.

\section{The Langlands correspondence for split tori}\label{cas-tores}

Split tori are isomorphic to $\mathbb G_{m}^{r}$, 
so there is nothing more than in the case of  $\mathbb G_{m}=GL_{1}$ explained in the previous section. 
Nevertheless the isomorphism of a split torus with  $\mathbb G_{m}^{r}$ is not canonical (because the automorphism group of $\mathbb G_{m}^{r}$
is non trivial, equal to $GL_{r}(\mathbb Z)$). Let $T$  be a split torus over $F$.  To obtain a canonical correspondence we introduce   the Langlands dual group $\widehat T$, defined as  the split torus over ${\mathbb Q}_\ell$ whose weights are the coweights of $T$ and reciprocally. In other words the lattice $\Lambda=\on{Hom}( \widehat T, \mathbb G_{m})$ is equal to  
$\on{Hom}(\mathbb G_{m},  T)$. 
Then the Langlands correspondence gives a bijection $\chi\mapsto \sigma(\chi)$ between 
\begin{itemize}
\item characters $\on{Bun}_{T}({{\mathbb F}_q})\to \overline{{\mathbb Q}_\ell}^{\times}$ with finite image
\item continuous morphisms $\pi_{1}(X,\overline \eta)\to \widehat T(\overline{{\mathbb Q}_\ell})$ with finite image
\end{itemize}
characterized by the fact that 
for any place $v$ of $X$ and any $\lambda \in \Lambda$ the image of 
$\sigma(\chi)(\on{Frob}_{v})$ by $\widehat T(\overline{{\mathbb Q}_\ell}) \xrightarrow{\lambda} \overline{{\mathbb Q}_\ell}^{\times}$  
 is equal to the image of $\mathcal O(v)$ by $\on{Pic}({{\mathbb F}_q})\xrightarrow{\lambda}
 \on{Bun}_{T}({{\mathbb F}_q})\xrightarrow{\chi} \overline{{\mathbb Q}_\ell}^{\times}$ (this condition is the particular case for tori of the condition of   
 ``compatibility with the Satake isomorphism'' which 
 we will consider later for all reductive groups). 
 
The construction of $\sigma(\chi)$ works as in the previous section, except that $a^{\mathbb Z}$ has to be replaced by a subgroup $\Xi$   of $\on{Bun}_{T}({{\mathbb F}_q})$ of finite index which is included in the kernel of $\chi$, and
we now have to use schemes of $T$-shtukas, defined 
using $T$-bundles instead of  line bundles.  

\section{Reminder on the dual group}  

Let $G$ be a split reductive group over $F$. 
            We denote by   $\widehat G$ the  Langlands dual group of   $G$. It is the split reductive group over   ${\mathbb Q}_\ell$ characterized by the fact that its roots and weights  are the coroots and coweights of    $G$, and reciprocally. 
  Here are some examples:  
         $$   \begin{array}{c | c}
G & \widehat G \\
\hline
GL_n & GL_{n} \\
        SL_n & PGL_n  \\
        SO_{2n+1} & Sp_{2n} \\
        Sp_{2n} & SO_{2n+1} \\
        SO_{2n} & SO_{2n}
        \end{array} $$
            and if  $G$ is one of the five exceptional groups, $\widehat G$ is of the same  type. Also the dual of a product of groups is the product of the dual groups.

 \begin{defi}     A global Langlands parameter is a conjugacy class of morphisms 
            $\sigma:\on{Gal}(\overline F/F)\to \widehat G(\overline{{\mathbb Q}_\ell})$ factorizing through 
            $\pi_{1}(U, \overline \eta)$  for some open dense $U\subset X$, 
            defined over a finite extension of ${\mathbb Q}_\ell$,  continuous and  semisimple. 
             \end{defi}                   
We say that $\sigma$ is semisimple if for any parabolic subgroup containing its image  there exists an associated  Levi subgroup containing it. 
Since $\overline{{\mathbb Q}_\ell}$ has characteristic $0$ this means equivalently that the Zariski  closure of  its image is reductive \cite{bki-serre}.

We now define the 
 Hecke operators (the spherical ones, also called unramified,  i.e.  without level). 
 They are similar to the Laplace operators on graphs. 

Let $v$ be a place of $X$. 
  If $\mathcal G$ and $\mathcal G'$ are two $G$-bundles over $X$ 
  we say that 
$(\mathcal G', \phi)$ is a  modification of $\mathcal G$ at $v$ if $\phi$ is an  
isomorphism between the  restrictions of $\mathcal G$ and $\mathcal G'$ to $X\setminus v$. 
Then the relative position   is a dominant coweight $\lambda$ of $G$ (in the case where $G=GL_{r}$ it is the  $r$-uple of elementary divisors). 
Let $\lambda$  be a  dominant coweight  of $G$.  We get the  Hecke correspondence
 $$   \xymatrix{ & \mathcal H_{v,\lambda}
         \ar[dl]^{h^{\leftarrow}} \ar[dr]_{h^{\rightarrow}} 
       &       \\
    \on{Bun}_{G}({{\mathbb F}_q}) &    & \on{Bun}_{G}({{\mathbb F}_q})    }$$     
where $\mathcal H_{v,\lambda}$ is the groupoid classifying 
modifications $(\mathcal G, \mathcal G',  \phi)$ at $v$ with relative position $\lambda$ 
and $h^{\leftarrow}$ and $h^{\rightarrow}$ send this object to $\mathcal G'$ and  $\mathcal G$. 
Then  the Hecke operator acts on functions by pullback by $h^{\leftarrow}$ followed by pushforward (i.e. sums in the fibers) by $h^{\rightarrow}$. In other words 
\begin{align*}T_{\lambda,v}: C_{c}^{\mathrm{cusp}}(\on{Bun}_{G}({{\mathbb F}_q})/\Xi,{\mathbb Q}_\ell) &  \to   C_{c}^{\mathrm{cusp}}(\on{Bun}_{G}({{\mathbb F}_q})/\Xi,{\mathbb Q}_\ell)   
 \\  f & \mapsto    \big[
\mathcal G\mapsto \sum_{(\mathcal G', \phi)
} f(\mathcal G')\big] \end{align*} 
where the finite sum is taken over  all the modifications $(\mathcal G', \phi)$ of $\mathcal G$ at $v$ 
with  relative position $\lambda$.   

These operators form an abstract commutative algebra $\mathcal H_{v}$, the so-called spherical (or unramified) Hecke  algebra at $v$, and this algebra acts on  $C_{c}^{\mathrm{cusp}}(\on{Bun}_{G}({{\mathbb F}_q})/\Xi,{\mathbb Q}_\ell)$. This algebra $\mathcal H_{v}$ is equal to  
$C_{c}(G(\mathcal O_{v})\backslash G(F_{v})/G(\mathcal O_{v}), {\mathbb Q}_\ell)$ and it is possible to write its action with the help of adèles. The actions of these  algebras $\mathcal H_{v}$ for different  $v$ commute with each other. 

The Satake isomorphism   \cite{satake, cartier-satake} can be viewed \cite{gross} as  a canonical   isomorphism 
$$[V]\mapsto T_{V,v}$$ from the Grothendieck ring of representations of  $\widehat G$ (with coefficients in $\mathbb Q_{\ell}$) to the unramifed Hecke algebra  $\mathcal H_{v}$, namely we have 
$T_{V\oplus V',v}=T_{V,v}+T_{V',v}$ and $T_{V\otimes V',v}=T_{V,v}T_{V',v}$. 
If $V$ is an irreducible representation of $\widehat G$, $T_{V,v}$ is a combination of the  $T_{\lambda,v}$ for $\lambda$ a weight of  $V$.  

\section{Presentation of the main result of \cite{coh}}\label{section-presentation}

 We now explain the construction in \cite{coh} of  a canonical decomposition  of  $\overline{{\mathbb Q}_\ell}$-vector spaces  
 \begin{gather}\label{dec-Cc-cusp}
C_{c}^{\mathrm{cusp}}(\on{Bun}_{G}({{\mathbb F}_q})/\Xi,\overline{{\mathbb Q}_\ell})=\bigoplus_{\sigma} \mathfrak H_{\sigma} \end{gather} 
where the direct sum is taken over  global Langlands parameters  $\sigma:\pi_{1}(X, \overline\eta)\to \widehat G(\overline{{\mathbb Q}_\ell})$. 
This  decomposition is  respected by and compatible with the action of  Hecke operators. 
In fact we  construct  a commutative ${\mathbb Q}_\ell$-algebra  
$$\mathcal B\subset \on{End}(C_{c}^{\mathrm{cusp}}(\on{Bun}_{G}({{\mathbb F}_q})/\Xi,{\mathbb Q}_\ell)) $$ 
containing  the image of $\mathcal H_{v}$ for all  places $v$ and such that 
  each character  $\nu$ of  $\mathcal B$ with values in $\overline{{\mathbb Q}_\ell}$ 
corresponds in a unique way to a    global Langlands parameter $\sigma$.

Since $\mathcal B$ is commutative we deduce  a   canonical  spectral decomposition 
   $$C_{c}^{\mathrm{cusp}}(\on{Bun}_{G}({{\mathbb F}_q})/\Xi,\overline{{\mathbb Q}_\ell})=\bigoplus_{\nu} \mathfrak H_{\nu} $$   where the direct sum is taken over characters $\nu$ of  $\mathcal B$ with values in  $\overline{{\mathbb Q}_\ell}$ and $\mathfrak H_{\nu} $ is the generalized eigenspace associated to $\nu$. 
  By associating to each $\nu$ a  global Langlands parameter $\sigma$  we deduce the  decomposition \eqref{dec-Cc-cusp} we want to construct. 
 We show in \cite{coh} that any $\sigma$ obtained in this way   factorizes through $\pi_{1}(X,\overline\eta)$, and that   the  decomposition \eqref{dec-Cc-cusp} is  compatible 
  with the Satake isomorphism at every place $v$ of $X$, in the sense that 
  for every representation $V$ of  $\widehat G$, $T_{V,v}$  acts on   
   $\mathfrak H_{\sigma} $ by multiplication by      $\on{Tr}_{V}(\sigma(\on{Frob}_{v}))$.

The elements of  $\mathcal B$ are constructed with the help of the $\ell$-adic cohomology  of  stacks of shtukas  and are called {\it excursion operators}. 

In the case of  $GL_{r}$,  since every semisimple linear representation is determined    up to conjugation by its character and since  the Frobenius elements $\on{Frob}_{v}$ are dense in 
$\on{Gal}(\overline F/F)$ by  the Chebotarev theorem, the  decomposition  \eqref{dec-Cc-cusp} is uniquely determined by its compatibility with the Satake isomorphism.

On the contrary, for some groups   $G$ other than $GL_{r}$, according to  Blasius and  Lapid \cite{blasius,lapid} it may happen that different   global Langlands parameters  correspond to the same characters of  $\mathcal H_{v}$ for every  place $v$. 
This comes from the fact that it is possible to find finite groups  $\Gamma$ 
and couples of   morphisms $\sigma,\sigma':\Gamma\to \widehat G(\overline{{\mathbb Q}_\ell})$ such that $\sigma$ and $\sigma'$ are not conjugated but that for any $\gamma\in \Gamma$, 
$\sigma(\gamma)$ and $\sigma'(\gamma)$ are  conjugated \cite{larsen1, larsen2}. 

Thus for a general group $G$, 
the algebra $\mathcal B$ of excursion operators may not be generated by the Hecke algebras $\mathcal H_{v}$ for all places $v$ and the 
compatibility of the decomposition  \eqref{dec-Cc-cusp} with Hecke operators may not characterize  it in a unique way. 
Therefore we wait for the construction  of the excursion operators (done in section \ref{section-main}) before we write the precise statement of our main result,  which will be theorem~\ref{intro-thm-ppal}.

  \section{The stacks of shtukas  and their $\ell$-adic cohomology}

The $\ell$-adic cohomology of a variety (over any algebraically closed field of characteristic $\neq \ell$) is very similar to the  Betti  cohomology of a complex variety, but it has coefficients in  ${\mathbb Q}_\ell$ (instead of  ${\mathbb Q}$ for the Betti cohomology). For its definition 
 Grothendieck introduced the notions of site and topos, which provide an extraordinary generalization of the usual notions of topological space and sheaf of sets on it. 

To a topological space $X$ we can associate the category  whose 
\begin{itemize}
\item objects are the open subsets $U\subset X$
\item    arrows $U\to V$ are the  inclusions $U\subset V$
\end{itemize}  and we have the    notion of a covering of  an open subset 
 by a family of open subsets. 
A site is an abstract category with a notion of covering 
of an object by a family of arrows targetting to it, 
with some natural axioms. 
A topos is the category of sheaves of sets on a site
(a sheaf of sets $\mathcal F$ on a site is a contravariant functor of ``sections of $\mathcal F$'' from the category 
of the site to the category of sets, satisfying, for each covering, a gluing axiom). 
Different sites may give the same topos.

To define the étale cohomology  of an algebraic variety $X$  we consider the étale  site 
\begin{itemize}
\item whose  objects are the  étale morphisms   $$\xymatrix{U
     \ar[d] 
     \\ X        }
$$   
\item  whose arrows are given by commutative triangles of étale morphisms,
$$\xymatrix{U \ar[dr]  \ar[rr]  &&V  \ar[dl] 
           \\  & X        }$$  
\item with the obvious  notion of covering.
\end{itemize}
  The  étale cohomology is defined with  cefficients in $\mathbb Z/\ell^{n}\mathbb Z$, whence $\mathbb Z_{\ell}$ by passing to the limit, and ${\mathbb Q}_{\ell}$ by inverting $\ell$.

The  stacks of shtukas, introduced by Drinfeld,  play a role analoguous to Shimura varieties over number fields. But they exist in a much greater generality. Indeed, while  the Shimura varieties are defined   over the spectrum of the ring of integers of a number field and are associated to  {\it a  minuscule coweight } of the dual group, the stacks of shtukas exist over {\it arbitrary powers} of the curve $X$, and can be associated to {\it arbitrary coweights}, as we will see now. One  simple reason for this difference between function fields and number fields is the following: in the  example of the product of two copies, the product $X\times X$ is taken over ${{\mathbb F}_q}$ whereas nobody knows what  the product 
$\on{Spec} \mathbb Z\times \on{Spec} \mathbb Z$ should be, and over what to take it.

Let $I$ be a finite set and $W=\boxtimes_{i\in I} W_{i}$ be an  irreducible ${\mathbb Q}_\ell$-linear representation of  $\widehat G^{I}$ (in other words each $W_{i}$ is an  irreducible representation of $\widehat G$).

We define $\on{Cht}_{I,W}$ as the reduced Deligne-Mumford  stack over  $X^{I}$ 
  whose points   over a scheme  $S$ over ${{\mathbb F}_q}$ classify shtukas, i.e.   
  \begin{itemize}
  \item   points  $(x_{i})_{i\in I}: S\to X^{I}$,  called the legs of the 
  shtuka (``les pattes du chtouca'' in French),  
\item  a    $G$-bundle $\mathcal G$ over  $ X\times S$,  
\item  an isomorphism 
$$\phi :\restr{\mathcal G }{(X\times S)\smallsetminus(\bigcup_{i\in I }\Gamma_{x_i})}\overset {\thicksim}{\to} \restr{(\on{Id}_{X}\times \on{Frob}_{S})^{*}(\mathcal G) }{(X\times S)\smallsetminus(\bigcup_{i\in I }\Gamma_{x_i})}$$ 
 (where   $\Gamma_{x_i}$ denotes the graph of $x_{i}$), such that 
 \begin{gather}\label{cond-position-relative}
\text{the relative position  at  } x_{i} \text{ of the  modification } \phi \text{ is bounded}
\\ \nonumber \text{by the  dominant coweight of   } G 
\text{ corresponding to the dominant weight of   } W_{i}.  
\end{gather}
  \end{itemize}
The notion of  Deligne-Mumford stack is in algebraic geometry what corresponds to the topological notion of orbifold. Every quotient of an algebraic variety by a finite étale group scheme 
is  a Deligne-Mumford stack and in fact $\on{Cht}_{I,W}$ is locally of this form.   
  \begin{rem}
 Compared to the notion of Artin stacks mentioned in remark~\ref{rem-artin},  a  Deligne-Mumford  stack is a particular case  where 
 the automorphism groups of geometric points are finite groups
 (instead of  algebraic groups).  
  \end{rem}

  \begin{rem} In the case of $GL_{1}$, {\it resp.} split tori, we had defined schemes  of shtukas.  With the above definition, the stacks of shtukas are the quotients 
  of these schemes by the trivial action of ${{\mathbb F}_q}^{\times}$, {\it resp.}   $T({{\mathbb F}_q})$. 
  \end{rem}

We denote by  $H_{I,W}$ the $\mathbb Q_\ell$-vector space equal to the ``Hecke-finite'' subspace of  the $\ell$-adic intersection  cohomology  
with compact support, in middle degree,  of the fiber  of  $\on{Cht}_{I,W}/\Xi$ over a generic geometric point of  $X^{I}$  (or,  in fact equivalently, over  a generic geometric point   of the  diagonal  $X\subset X^{I}$).   To give an idea of intersection cohomology, let us say that for a smooth variety it is the same as the $\ell$-adic cohomology and that for (possibly singular) projective varieties it is  Poincaré self-dual. 
An element of this $\ell$-adic  intersection  cohomology is said to be Hecke-finite if it belongs to a sub-$\mathbb Z_{\ell}$-module of finite type  stable by all Hecke operators $T_{\lambda,v}$ (or equivalently by all Hecke operators $T_{V,v}$). Hecke-finiteness is a technical condition  but Cong Xue has proven \cite{these-cong} that $H_{I,W}$ can equivalently be defined by a cuspidality condition (defined using stacks of shtukas for parabolic subgroups of $G$ and their Levi quotients) and that it has finite dimension over $\mathbb Q_\ell$.

Drinfeld has constructed  ``partial Frobenius morphisms'' between stacks of shtukas. To define them we need  
 a small generalization of the stacks  of shtukas 
  where we require a factorization of  $\phi$ as a composition of several  modifications. Let $(I_{1},...,I_{k})$ be  an ordered partition    of $I$. 
 An example is the coarse partition $(I)$ and in fact the stack $\on{Cht}_{I,W} $ previously defined is equal to  $\on{Cht}_{I,W}^{(I)}$ in the following definition. 
 
 \begin{defi}\label{defi-Cht-I1-Ik}
 We define      $\on{Cht}_{I,W} ^{(I_{1},...,I_{k})}$ as the {\it reduced }   Deligne-Mumford  stack    whose points over 
 a scheme   $S$ over  ${{\mathbb F}_q}$  classify   \begin{gather}\label{intro-donnee-chtouca}\big( (x_i)_{i\in I}, \mathcal G_{0} \xrightarrow{\phi_{1}}  \mathcal G_{1}  \xrightarrow{\phi_{2}}
\cdots\xrightarrow{\phi_{k-1}}   \mathcal G_{k-1}  \xrightarrow{ \phi_{k}}    (\on{Id}_{X}\times \on{Frob}_{S})^{*}(\mathcal G_{0}) 
\big)
\end{gather}
with 
 \begin{itemize}
\item $x_i\in (X\smallsetminus N)(S)$ for $i\in I$, 
\item for $i\in \{0,...,k-1\}$, $\mathcal G_{i}$   is a   $G$-bundle over  $X\times S$ and we write 
   $\mathcal G_{k} =(\on{Id}_{X}\times \on{Frob}_{S})^{*}( \mathcal G_{0})$ to prepare the next item, 
 \item  
for   $j\in\{1,...,k\}$
 $$\phi_{j}:\restr{\mathcal G_{j-1}}{(X\times S)\smallsetminus(\bigcup_{i\in I_{j}}\Gamma_{x_i})}\overset {\thicksim}{\to} \restr{\mathcal G_{j}}{(X\times S)\smallsetminus(\bigcup_{i\in I_{j}}\Gamma_{x_i})}$$ is an   isomorphism  such that the  relative position  of  $\mathcal G_{j-1}$ with respect to $\mathcal G_{j}$ at  $x_{i}$ (for  $i\in I_{j}$) is bounded by the   dominant coweight  of  $G$ corresponding to the  dominant weight  of  $W_{i}$. 
 \end{itemize}
 \end{defi}

We can show that the obvious  morphism  $\on{Cht}_{I,W} ^{(I_{1},...,I_{k})}\to \on{Cht}_{I,W}$ (which forgets the intermediate modifications $\mathcal G_{1}, ..., \mathcal G_{k-1}$) gives an isomorphism at the level of intersection cohomology. The interest of $\on{Cht}_{I,W} ^{(I_{1},...,I_{k})}$ is that we have the partial Frobenius morphism 
$\on{Frob}_{I_{1}}: \on{Cht}_{I,W} ^{(I_{1},...,I_{k})} \to \on{Cht}_{I,W} ^{(I_{2},...,I_{k}, I_{1})} $ which sends  \eqref{intro-donnee-chtouca} to 
 \begin{gather}\nonumber
 \big( (x'_i)_{i\in I},   \mathcal G_{1}  \xrightarrow{\phi_{2}}
\cdots\xrightarrow{\phi_{k-1}}   \mathcal G_{k-1}  \xrightarrow{ \phi_{k}}    (\on{Id}_{X}\times \on{Frob}_{S})^{*}(\mathcal G_{0}) \\ \nonumber 
\xrightarrow{ (\on{Id}_{X}\times \on{Frob}_{S})^{*}(\phi_{1})}    (\on{Id}_{X}\times \on{Frob}_{S})^{*}(\mathcal G_{1}) \big)
\end{gather}
where $x_{i}'=\on{Frob}(x_{i})$ if $i\in I_{1}$ and $x_{i}'=x_{i}$ otherwise. 
Taking $I_{1}$ to be a singleton we get the action on $H_{I,W}$  of the 
partial Frobenius morphisms. 
Thanks to an extra work (using the  Hecke-finiteness condition and Eichler-Shimura relations), we are able in \cite{coh} to apply Drinfeld's lemma, and this endows the  $\mathbb Q_\ell$-vector space $H_{I,W}$   with a continuous action of $\on{Gal}(\overline F/F)^{I}$. 
  
For  $I=\emptyset$ and  $W=\mathbf  1$ (the trivial representation), we have   \begin{gather}\label{iso-empty}       H_{\emptyset,\mathbf  1}=C_{c}^{\mathrm{cusp}}(\on{Bun}_{G}({{\mathbb F}_q})/\Xi,\mathbb Q_\ell). \end{gather}
     Indeed the $S$-points over $\on{Cht}_{\emptyset,\mathbf  1}$ classify the  
 $G$-bundles $\mathcal G$ over  $ X\times S$, 
equipped with an  isomorphism 
$$\phi : \mathcal G  \overset {\thicksim}{\to}  (\on{Id}_{X}\times \on{Frob}_{S})^{*}(\mathcal G)  . $$ 
If we see  $\mathcal G$ as a  $S$-point of  $\on{Bun}_{G}$, 
$(\on{Id}_{X}\times \on{Frob}_{S})^{*}(\mathcal G)$ is its image by 
$\on{Frob}_{\on{Bun}_{G}}$. Therefore  $\on{Cht}_{\emptyset,\mathbf  1}$
classifies the fixed points of  $\on{Frob}_{\on{Bun}_{G}}$ and it is dicrete (i.e. of dimension $0$) and  equal to 
$ \on{Bun}_{G} ({{\mathbb F}_q})$. Therefore the $\ell$-adic cohomology of $\on{Cht}_{\emptyset,\mathbf  1}/\Xi$ is equal to 
$C_{c}(\on{Bun}_{G}({{\mathbb F}_q})/\Xi,\mathbb Q_\ell)$
and  in this particular case it is easy to see  that  Hecke-finiteness  is equivalent to cuspidality, so that 
\eqref{iso-empty}  holds true. 
 
Up to now  we defined a  vector space $H_{I,W}$ for every  isomorphism class of  irreducible representation $W=\boxtimes_{i\in I} W_{i}$  of  $\widehat G^{I}$.
A construction based on the geometric Satake equivalence   enables to   \begin{itemize}
  \item [] {\bf a)} define $H_{I,W}$ functorialy in $W$
  \item [] {\bf b)}  understand the fusion of legs   \end{itemize}
as  explained in the next  proposition.

   \begin{prop} \label{prop-abc}
     a) For every finite set   $I$,      $$W\mapsto  
    H_{I,W},  \ \ u\mapsto \mathcal H(u)$$  is a  $\mathbb Q_\ell$-linear  functor from the category of finite dimensional representations of   $\widehat G^{I}$ to the category of finite dimensional and continuous representations of  $\on{Gal}(\overline F/F)^{I}$.      
    
    This means that for every morphism 
    $$u:W\to W'$$ of representations of  $\widehat G^{I}$, we have a  morphism 
    $$ \mathcal H(u): H_{I,W}\to H_{I,W'}$$ of  representations   of      $\on{Gal}(\overline F/F)^{I}$.

b) 
For each map   $\zeta: I\to J$ between finite sets, 
 we have an  isomorphism 
    $$
           \chi_{\zeta}: H_{I,W}\overset {\thicksim}{\to} 
 H_{J,W^{\zeta}}  $$
 which is 
 \begin{itemize}
 \item  functorial in   $W$, where  $W$ is a   representation of $\widehat G^{I}$ and   $W^{\zeta}$ denotes the  representation of $\widehat G^{J}$ on  $W$ obtained by composition with the diagonal   morphism  $$  \widehat G^{J}\to \widehat G^{I}, (g_{j})_{j\in J}\mapsto (g_{\zeta(i)})_{i\in I} $$ 
 \item $\on{Gal}(\overline F/F)^{J}$-equivariant, where $\on{Gal}(\overline F/F)^{J}$ acts on the LHS by the  diagonal  morphism  
$$\on{Gal}(\overline F/F)^{J}\to \on{Gal}(\overline F/F)^{I},  \ (\gamma_{j})_{j\in J}\mapsto (\gamma_{\zeta(i)})_{i\in I}, 
$$
 \item    and compatible with  composition, i.e.  for every 
 $I\xrightarrow{\zeta} J\xrightarrow{\eta} K$ we have 
 $\chi_{\eta\circ \zeta}=\chi_{\eta}\circ\chi_{\zeta}$. 
   \end{itemize}
           \end{prop}

The statement  b) is a bit complicated, here is a basic example of it. 
 For every finite set  $I$ we write $\zeta_{I}:I\to \{0\}$ the tautological map 
(where $\{0\}$ is an arbitrary choice of  notation for a singleton). 
 If $W_{1}$ and $W_{2}$ are two representations of $\widehat G$, 
the statement of  b) provides a  canonical isomorphism 
\begin{gather} \label{isom-fusion}\chi_{\zeta_{\{1,2\}}}:H_{\{1,2\},W_{1}\boxtimes W_{2}}\overset {\thicksim}{\to} H_{\{0\},W_{1}\otimes W_{2}}\end{gather}
 associated to $\zeta_{ \{1,2\} }: \{1,2\} \to \{0\} $. 
 We stress the difference  between  $W_{1}\boxtimes W_{2}$ which is a  representation of  $(\widehat G)^{2}$ and $W_{1}\otimes W_{2}$ which is a   representation of $\widehat G$.

Another example of  b) is the isomorphism on the left in 
 \begin{gather} \label{iso-vide-singleton-cusp} H_{\{0\},\mathbf 1}\isor{\chi_{\zeta_{\emptyset}}^{-1}} H_{ \emptyset ,\mathbf 1}
 \overset{\eqref{iso-empty}}{=} C_{c}^{\mathrm{cusp}}(\on{Bun}_{G}({{\mathbb F}_q})/\Xi,\mathbb Q_\ell) \end{gather}    
  which is associated to  $\zeta_{\emptyset}: \emptyset \to \{0\}$   (the idea of the isomorphism $\chi_{\zeta_{\emptyset}}$ is that $H_{ \emptyset ,\mathbf 1}$ 
{\it resp.} $H_{\{0\},\mathbf 1}$ is the 
 cohomology of the stack of  shtukas without legs,  {\it resp.} with a inactive leg,  and that they are equal).  
Thanks to  \eqref{iso-vide-singleton-cusp}
  we are reduced to construct a  decomposition 
    \begin{gather}\label{dec-singleton} H_{\{0\},\mathbf 1}\otimes_{{\mathbb Q}_\ell} \overline{{\mathbb Q}_\ell}=\bigoplus_{\sigma} \mathfrak H_{\sigma}. \end{gather}

{\noindent \bf Idea of the proof of proposition \ref{prop-abc}. } We denote by   $\on{Cht}_{I}$   the inductive limit of   Deligne-Mumford stacks over  $X^{I}$,  defined as 
  $\on{Cht}_{I,W}$ above, but without the  condition \eqref{cond-position-relative} on the  relative position. In other words, 
  and with an extra letter $\mathcal G'$ to prepare the next definition, 
  the  points  of $\on{Cht}_{I}$ over a scheme  $S$ over ${{\mathbb F}_q}$ classify   
  \begin{itemize}
  \item   points  $(x_{i})_{i\in I}: S\to X^{I}$,   
  \item  two   $G$-bundles $\mathcal G$ and $\mathcal G'$ over  $ X\times S$,  
\item   a modification $\phi$ at the  $x_{i}$, i.e.  an isomorphism 
$$\phi :\restr{\mathcal G }{(X\times S)\smallsetminus(\bigcup_{i\in I }\Gamma_{x_i})}\overset {\thicksim}{\to} \restr{\mathcal G' }{(X\times S)\smallsetminus(\bigcup_{i\in I }\Gamma_{x_i})}$$ 
  \item an isomorphism $\theta: \mathcal G' \overset {\thicksim}{\to} (\on{Id}_{X}\times \on{Frob}_{S})^{*}(\mathcal G)$. 
  \end{itemize}
We introduce the     ``prestack'' $\mathcal M_{I}$ of  ``modifications on the formal neighborhood of the  $x_{i}$'', whose points   over a scheme   $S$ over  ${{\mathbb F}_q}$ classify 
  \begin{itemize}
  \item   points  $(x_{i})_{i\in I}: S\to X^{I}$,   
  \item  two   $G$-bundles $\mathcal G$ and $\mathcal G'$ on the formal completion   $\widehat{ X\times S}$ of $X\times S$ in the  neighborhood  of the union of the graphs $\Gamma_{x_i}$,  
\item   a modification $\phi$ at the $x_{i}$, i.e.  an isomorphism 
$$\phi :\restr{\mathcal G }{(\widehat{ X\times S})\smallsetminus(\bigcup_{i\in I }\Gamma_{x_i})}\overset {\thicksim}{\to} \restr{\mathcal G' }{(\widehat{ X\times S})\smallsetminus(\bigcup_{i\in I }\Gamma_{x_i})}. $$    \end{itemize}
The expert  reader will notice that  for any morphism $S\to X^{I}$, $ \mathcal M_{I}\times _{X^{I}}S$ is the  quotient of the   affine grassmannian of Beilinson-Drinfeld over $S$ by $\Gamma(\widehat{ X\times S}, G)$. 
 We have a formally smooth  morphism  $\epsilon_{I}:\on{Cht}_{I}\to \mathcal M_{I}$ given by  restricting 
$\mathcal G$ and $\mathcal G'$ to the formal neighborhood of the graphs of the  $x_{i} $ and forgetting  $\theta$.

The geometric Satake equivalence,  
due to  Lusztig, Drinfeld, Ginzburg and  Mirkovic-Vilonen 
   \cite{hitchin,mv},  
 is a fundamental statement 
which {\it  constructs}  $\widehat G$ from  $G$ and is the cornerstone of the geometric 
  Langlands program. It is a  canonical equivalence of tensor categories between 
 \begin{itemize}
 \item  the category of perverse sheaves on the fiber of  $\mathcal M_{\{0\}}$ above any point of $X$ (where $\{0\}$ is an arbitrary notation for a singleton)
    \item the tensor category of  representations of $\widehat G$. 
 \end{itemize}
 
 For the non expert reader we recall that perverse sheaves, introduced in \cite{bbd}, behave like ordinary sheaves and have, in spite of their name, very good properties.  An example is given by   intersection cohomology sheaves of closed (possibly singular) subvarieties, whose total cohomology is the intersection cohomology of this subvarieties. 
 
   The tensor structure on the first category above is obtained by 
     ``fusion of legs'', thanks to the fact that  $\mathcal M_{\{1,2\}}$  is equal to  $\mathcal M_{\{0\}}\times \mathcal M_{\{0\}}$ outside the diagonal of  $X^{2}$ and to  $\mathcal M_{\{0\}}$ on the diagonal. 
     The first category is tannakian and $\widehat G$ is defined as the group of automorphisms of a natural fiber functor. 
     
    This equivalence gives,  for every  representation $W$ of $\widehat G^{I}$,  a perverse sheave $\mathcal S_{I,W}$ on   $\mathcal M_{I}$, with  the following properties: 
   \begin{itemize}
   \item $\mathcal S_{I,W}$  is functorial in  $W$, 
   \item for every  surjective  map $I\to J$, $\mathcal S_{J,W^{\zeta}}$ is canonically isomorphic to the  restriction of  $\mathcal S_{I,W}$ to  $\mathcal M_{I}\times_{X^{I}}X^{J}\simeq \mathcal M_{J}$, where  $X^{J}\to X^{I}$ is the  diagonal morphism, 
   \item for every   irreducible representation $W$, $\mathcal S_{I,W}$ is the intersection cohomology sheaf of the closed substack of  $\mathcal M_{I}$ defined by the condition \eqref{cond-position-relative} on the  relative position  of the  modification $\phi$ at the  $x_{i}$. 
   \end{itemize}
  Then we define  $H_{I,W}$ as the  ``Hecke-finite'' subspace of the   cohomology  
with compact support of  $\epsilon_{I}^{*}(\mathcal S_{I,W})$ 
  on the fiber of  $\on{Cht}_{I}/\Xi$ over a geometric generic point of  $X^{I}$  (or, in fact equivalently, over a geometric generic point of the diagonal $X\subset X^{I}$). 
  The first two properties above imply a) and b) of the proposition. The third one and the smoothness  of $\epsilon_{I}$ ensure that,  for  $W$ irreducible,
  $\epsilon_{I}^{*}(\mathcal S_{I,W})$ is the intersection cohomology sheaf of $\on{Cht}_{I,W}$ and therefore 
    the new definition of   $H_{I,W}$  generalizes the first one using  the intersection cohomology of  $\on{Cht}_{I,W}$.

\section{Excursion operators and the main theorem of \cite{coh}}\label{section-main}
  
  Let $I$ be a finite set. Let $(\gamma_{i})_{i\in I}\in \on{Gal}(\overline F/F)^{I}$. Let    $W$ be a  representation  of  $\widehat G^{I}$ and  $x\in W$ and $\xi\in W^{*}$ be invariant by the diagonal action of   $\widehat G$. 
      We define the endomorphism 
    $S_{I,W,x,\xi,(\gamma_{i})_{i\in I}}$ of  
 \eqref{iso-vide-singleton-cusp} 
as the composition  
  \begin{gather}\label{excursion-def-intro}
  H_{\{0\},\mathbf  1}\xrightarrow{\mathcal H(x)}
 H_{\{0\},W^{\zeta_{I}}}\isor{\chi_{\zeta_{I}}^{-1}} 
  H_{I,W}
  \xrightarrow{(\gamma_{i})_{i\in I}}
  H_{I,W} \isor{\chi_{\zeta_{I}}} H_{\{0\},W^{\zeta_{I}}}  
  \xrightarrow{\mathcal H(\xi)} 
  H_{\{0\},\mathbf  1} 
    \end{gather}
    where $\mathbf 1$ denotes the trivial  representation of $\widehat G$, 
and  $x:\mathbf 1\to W^{\zeta_{I}}$ and $\xi: W^{\zeta_{I}}\to \mathbf 1$ are considered as  morphisms of representations of $\widehat G$ (we recall that $\zeta_{I}:I\to \{0\}$ is the obvious map and that  $W^{\zeta_{I}}$ is simply the vector space  $W$ equipped with the diagonal  action of $\widehat G$).

Paraphrasing \eqref{excursion-def-intro} this  operator is the composition 
    \begin{itemize}
    \item of a creation   operator associated to   $x$, whose effect is to create legs at the same (generic) point of the curve, 
    \item  of a Galois action, 
    which moves the legs on the curve independently from each other, 
   then brings them back to the same (generic) point of the curve, 
       \item of an  annihilation  operator associated to  $\xi$. 
         \end{itemize}
It is called  an  ``excursion operator'' because it moves the legs on the curve  (this is what makes it  non trivial).

  To $W,x,\xi$ we associate the matrix  coefficient   $f$  defined by 
       \begin{gather} \label{rel-f-Wxxi}   f ((g_{i})_{i\in I})=\langle \xi, (g_{i})_{i\in I}\cdot x \rangle.\end{gather}  
We see that   $f$ is a  function on $\widehat G^{I}$ invariant  by left and right  translations by the diagonal $\widehat G$. In other words  $f\in\mathcal O(\widehat G\backslash \widehat G^{I}/\widehat G)$,  where $\widehat G\backslash \widehat G^{I}/\widehat G$ denotes the coarse quotient, defined as the spectrum of the algebra of functions   $f$ as above. 
Unlike the stacky quotients considered before, the coarse quotients are schemes  and therefore forget the automorphism groups of points.

For every  function $f\in\mathcal O(\widehat G\backslash \widehat G^{I}/\widehat G)$ we can find  $W,x,\xi$ such that  \eqref{rel-f-Wxxi}  holds. 
   We show easily that    $S_{I,W,x,\xi,(\gamma_{i})_{i\in I}}$  does not depend on the choice of  $W,x,\xi$ satisfying  \eqref{rel-f-Wxxi}, and therefore we denote it by  
    $S_{I,f,(\gamma_{i})_{i\in I}}$.

The  conjectures of Arthur and Kottwitz on multiplicities in vector spaces of  automorphic forms and in the  cohomologies of    Shimura varieties \cite{arthur-ast,kottwitz3} 
give,  by extrapolation to stacks of shtukas,  the following heuristics. 
 
  \begin{rem}  Heuristically we  conjecture 
    that for every  global Langlands parameter $\sigma$
there exists a  $\overline{{\mathbb Q}_\ell}$-linear representation $\mathfrak A_{\sigma}$ of its centralizer $S_{\sigma}\subset \widehat G$ (factorizing through $S_{\sigma}/Z_{\widehat G}$), so that  we have a $\on{Gal}(\overline F/F)^{I}$-equivariant isomorphism 
  \begin{gather}\label{heuri} H_{I,W}\otimes_{{\mathbb Q}_\ell} \overline{{\mathbb Q}_\ell} 
  \overset{?}{=}\bigoplus_{\sigma} \big( \mathfrak A_{\sigma}\otimes_{\overline{{\mathbb Q}_\ell}} W_{\sigma^{I}}\big)^{S_{\sigma}} \end{gather}
   where $W_{\sigma^{I}}$ is the $\overline{{\mathbb Q}_\ell}$-linear representation of  
 $\on{Gal}(\overline F/F)^{I}$ obtained by composition of the  representation $W$ of  
 $\widehat G^{I}$ with the  morphism $\sigma^{I}: \on{Gal}(\overline F/F)^{I}\to \widehat G(\overline{{\mathbb Q}_\ell})^{I}$, and $S_{\sigma}$ acts diagonally. 
  We conjecture that  \eqref{heuri}   is functorial in $W$, compatible to $\chi_{\zeta}$ and that it is equal to the decomposition \eqref{dec-singleton}   when  $W=\bf 1$ 
  (so that $\mathfrak H_{\sigma}=
  ( \mathfrak A_{\sigma})^{S_{\sigma}}$). 
 
In the heuristics  \eqref{heuri}  the endomorphism 
    $S_{I,f,(\gamma_{i})_{i\in I}}=S_{I,W,x,\xi,(\gamma_{i})_{i\in I}}$ of   
  $$H_{\{0\},\mathbf 1}\otimes_{{\mathbb Q}_\ell} \overline{{\mathbb Q}_\ell} 
  \isor{\chi_{\zeta_{\emptyset}}^{-1}} H_{ \emptyset ,\mathbf 1}\otimes_{{\mathbb Q}_\ell} \overline{{\mathbb Q}_\ell}
  \overset{?}{=}\bigoplus_{\sigma} ( \mathfrak A_{\sigma})^{S_{\sigma}}$$    
 acts on  $( \mathfrak A_{\sigma})^{S_{\sigma}}$ by the composition 
  \begin{gather*}
 ( \mathfrak A_{\sigma})^{S_{\sigma}}
 \isor{\chi_{\zeta_{\emptyset}}} 
 ( \mathfrak A_{\sigma} \otimes  \mathbf 1)^{S_{\sigma}}
 \xrightarrow{\on{Id}_{  \mathfrak A_{\sigma}}\otimes x}
 ( \mathfrak A_{\sigma}  \otimes W_{\sigma^{I}})^{S_{\sigma}}
  \xrightarrow{(\sigma(\gamma_{i}))_{i\in I} }
  ( \mathfrak A_{\sigma}  \otimes W_{\sigma^{I}} )^{S_{\sigma}}
  \\ \nonumber 
  \xrightarrow{\on{Id}_{ \mathfrak A_{\sigma}}\otimes \xi} 
 ( \mathfrak A_{\sigma}  \otimes  \mathbf 1)^{S_{\sigma}}
 \isor{\chi_{\zeta_{\emptyset}}^{-1}} 
  ( \mathfrak A_{\sigma})^{S_{\sigma}}
    \end{gather*}
  i.e.  by the scalar 
  $$ \langle{\xi, (\sigma(\gamma_{i}))_{i\in I}\cdot x} \rangle =f\big( (\sigma(\gamma_{i}))_{i\in I} \big).
   $$
  In other words we should have 
  \begin{gather}\label{heuri-Hsigma}
  \mathfrak H_{\sigma}\overset{?}{=}
  \text{eigenspace of the  }S_{I,f,(\gamma_{i})_{i\in I}}\text{ with the eigenvalues }f\big( (\sigma(\gamma_{i}))_{i\in I} \big).\end{gather} 
\end{rem}

The heuristics  \eqref{heuri-Hsigma}
clearly indicates the path to follow. 
We show in \cite{coh} that the  $S_{I,f,(\gamma_{i})_{i\in I}}$ generate a commutative ${\mathbb Q}_\ell$-algebra $\mathcal B$ and satisfy some properties implying
the following proposition. 

  \begin{prop} \label{prop-existence-sigma} For each character  $\nu$ of  $\mathcal B$ with values in  $\overline{{\mathbb Q}_\ell}$ there exists a  unique global Langlands parameter $\sigma$ such that for all  $I,f$ and   $(\gamma_{i}  )_{i\in I}$, we have 
  \begin{gather}\label{caract-sigma}\nu(S_{I,f,(\gamma_{i})_{i\in I}})= f ((\sigma(\gamma_{i}))_{i\in I}).
   \end{gather}
 \end{prop}
The unicity of $\sigma$ in the previous proposition comes from the fact that, for any integer $n$, taking $I=\{0,...,n\}$, the coarse quotient  
  $\widehat G\backslash \widehat G^{I}/\widehat G$ identifies with the coarse quotient of 
  $(\widehat G)^{n}$ by diagonal conjugation by $\widehat G$, and therefore, for any 
  $(\gamma_{1}, ..., \gamma_{n})\in (\on{Gal}(\overline F/F))^{n}$, 
  \eqref{caract-sigma} applied to $(\gamma_{i})_{i\in I}=(1,\gamma_{1}, ..., \gamma_{n})$ determines $(\sigma(\gamma_{1}), ..., \sigma(\gamma_{n}))$ up to conjugation and semisimplification (thanks to \cite{richardson}). 
  The existence and continuity of $\sigma$ are justified thanks to relations and topological properties satisfied by  the excursion operators.

   Since $\mathcal B$ is commutative we obtain a   canonical  spectral decomposition 
   $   H_{\{0\},\mathbf 1} \otimes_{{\mathbb Q}_\ell} \overline{{\mathbb Q}_\ell} =\bigoplus_{\nu} \mathfrak H_{\nu} $
   where the direct sum is taken over  characters $\nu$ of  $\mathcal B$ with values in  $\overline{{\mathbb Q}_\ell}$.  
  Associating to each  $\nu$ a unique global Langlands parameter $\sigma$ as in the previous proposition,  we deduce the  decomposition \eqref{dec-singleton} we wanted to construct. We do not know if $\mathcal B$ is reduced.

   Moreover   the  unramified Hecke operators are particular cases of excursion operators:  for every place $v$ and 
   for every irreducible representation  $V$ of  $\widehat G$ with character $\chi_{V}$, the  unramified Hecke operator $T_{V,v}$  
      is equal to the excursion operator $ S_{\{1,2\}, f,(\on{Frob}_{v},1)}$ where $f\in \mathcal O(\widehat G\backslash (\widehat G)^{2}/\widehat G)$ is given by  $f(g_{1},g_{2})=\chi_{V}(g_{1}g_{2}^{-1})$, and $\on{Frob}_{v}$ is  a Frobenius element at  $v$. 
          This    is proven in \cite{coh} by a geometric argument  (essentially a computation of the intersection of algebraic cycles in a stack of  shtukas). 
It implies the compatibility of the decomposition \eqref{dec-singleton} with the  Satake isomorphism at all  places.  

\begin{rem}
By the Chebotarev density theorem, the subalgebra of $\mathcal B$ generated by all the Hecke algebras $\mathcal H_{v}$ is equal to the subalgebra generated by 
the excursion operators with $\sharp I=2$. The remarks at the end of section \ref{section-presentation} show that in general  it is necessary to consider excursion operators with $\sharp I>2$ to generate the whole algebra $\mathcal B$. 
\end{rem}

Finally we can state the main theorem.  

\begin{thm}  \label{intro-thm-ppal}  We have a canonical    decomposition 
of  $\overline{{\mathbb Q}_\ell}$-vector spaces 
 \begin{gather} \label{dec-intro-thm-ppal}
 C_{c}^{\mathrm{cusp}}(\on{Bun}_{G}({{\mathbb F}_q})/\Xi,\overline{{\mathbb Q}_\ell})=\bigoplus_{\sigma}
 \mathfrak H_{\sigma},\end{gather}
 where the direct sum in the RHS is indexed by   global Langlands parameters , i.e.   $\widehat G(\overline{{\mathbb Q}_\ell})$-conjugacy  classes of  morphisms 
       $\sigma:\on{Gal}(\overline F/F)\to \widehat G(\overline{{\mathbb Q}_\ell})$   factorizing through $\pi_{1}(X,\overline\eta)$, 
      defined over a finite extension of  ${\mathbb Q}_\ell$, continuous and   semisimple. 
          
   This  decomposition  is uniquely determined by the following property : 
       $ \mathfrak H_{\sigma}$ is equal to the generalized eigenspace   
       associated to the character  $\nu$ of $\mathcal B$ defined by  
       \begin{gather}\label{relation-fonda}\nu(S_{I,f,(\gamma_{i})_{i\in I}})=f((\sigma(\gamma_{i}))_{i\in I}. \end{gather}
       
     This decomposition is respected by the  Hecke operators and is  compatible with the Satake isomorphism at all  places  $v$  of  $X$.     \end{thm}

Everything is still true with a level  (a finite subscheme  $N$ of $X$). 
We denote by $\on{Bun}_{G,N}({{\mathbb F}_q})$ the set of isomorphism classes  of  $G$-bundles over  $X$ trivialized on  $N$. 
Then we have a  canonical decomposition 
 \begin{gather} \label{dec-nivau-N}
 C_{c}^{\mathrm{cusp}}(\on{Bun}_{G,N}({{\mathbb F}_q})/\Xi,\overline{{\mathbb Q}_\ell})=\bigoplus_{\sigma}
 \mathfrak H_{\sigma},\end{gather}
 where the direct sum  is taken over   global Langlands parameters   $\sigma:\pi_{1}(X\smallsetminus N,\overline\eta)\to \widehat G(\overline{{\mathbb Q}_\ell})$. This  decomposition is 
 respected by all Hecke operators and compatible with the  Satake isomorphism at all  places of $X\smallsetminus N$.
If  $G$ is split we have, by \cite{thang},    \begin{gather}\label{egalite-deploye}\on{Bun}_{G,N}({{\mathbb F}_q})=G(F)\backslash G(\mathbb A)/K_{N}\end{gather} 
(where $\mathbb A$ is the ring of  adèles, 
$\mathbb O$ is the ring of  integral adèles, $\mathcal O_{N}$ the ring of  functions  on  $N$ and 
$K_{N}=\on{Ker}(G(\mathbb O)\to G(\mathcal O_{N}))$). 
 When  $G$ is non necessarily split the RHS of  \eqref{egalite-deploye} must be replaced by a direct sum, indexed by the finite group  $ \on{ker}^{1}(F,G)$,  
 of adelic quotients for inner forms of  $G$ and in the definition of   global Langlands parameters  
we must replace  $\widehat G$ by the $L$-group (see \cite{borel-corvallis} for $L$-groups). 
 
 We have a statement similar to theorem~\ref{intro-thm-ppal} with coefficients in $\overline{{\mathbb F}_\ell}$ instead of $\overline{{\mathbb Q}_\ell}$. 
 
We can also consider the case of  metaplectic groups thanks to the 
 metaplectic variant of the geometric  Satake equivalence due to  Finkelberg and Lysenko \cite{finkelberg-lysenko,lysenko-red,dennis-sergey}. 

\begin{rem}\label{rem-drinfeld-construction}
Drinfeld gave an idea to prove something like the heuristics  \eqref{heuri} but it is a bit difficult to formulate the result. 
Let  $\on{Reg}$ be the left regular  representation  of $\widehat G$ with coefficients in ${\mathbb Q}_\ell$ (considered as an inductive limit of finite dimensional representations). 
We can endow  the ${\mathbb Q}_\ell$-vector space $H_{\{0\},\on{Reg}}$ (of infinite dimension in general) with 
 \begin{itemize}
 \item  [a) ] 
 a structure of  
  module over the algebra of functions on the  ``affine space $\mathcal S$ of  morphisms $\sigma:\on{Gal}(\overline F/F)\to 
 \widehat G$ with coefficients in ${\mathbb Q}_\ell$-algebras'', 
 \item [b)]  an  algebraic action of $\widehat G$ (coming from the right action of 
 $\widehat G$ on $\on{Reg}$) which is 
 compatible with conjugation by  $\widehat G$ 
on $\mathcal S$. 
 \end{itemize}
 The space $\mathcal S$ is not rigorously defined and the rigorous definition of  structure a) is 
 the following. For any finite dimensional $\mathbb Q_{\ell}$-linear representation $V$ of $\widehat G$, with underlying vector space $\underline V$, 
  $H_{\{0\},\on{Reg}}\otimes \underline V$ is equipped with an action of $\on{Gal}(\overline F/F)$, making it an inductive limit of finite dimensional continuous representations of $\on{Gal}(\overline  F/F)$, as follows. 
We have a  $\widehat  G$-equivariant isomorphism 
   \begin{align*} \theta:  \on{Reg} \otimes \underline V  & \simeq  \on{Reg}\otimes V  \\
     f\otimes x & \mapsto    [g\mapsto f(g)  g.x]
  \end{align*} 
    where  $\widehat  G$ acts diagonally on the RHS, and to give a meaning to the formula the RHS is identified with the vector space of algebraic functions $\widehat G\to V$. 
  Therefore we have an  isomorphism 
    $$  H_{\{0\}, \on{Reg} } \otimes  \underline   V = H_{\{0\}, \on{Reg}  \otimes \underline V }
    \isor{\theta} H_{\{0\},\on{Reg}  \otimes V}
    \simeq H_{\{0,1\},\on{Reg} \boxtimes V}$$
    where the first equality is tautological (since $ \underline V$ is just a vector space) and the last isomorphism is the inverse of the fusion isomorphism $\chi_{\zeta_{\{0,1\}}} $ of \eqref{isom-fusion}. 
     Then the  action of  $ \on{Gal}(\overline F/F)$ on  
the LHS   is defined as  the action  of  $ \on{Gal}(\overline F/F)$ on the RHS  corresponding to the leg $1$.  
 If   $V_{1}$ and  $V_{2}$ are two representations of  $\widehat G$, 
 the two actions of  $ \on{Gal}(\overline F/F)$ on  $H_{\{0\},\on{Reg}}\otimes  \underline  {V_{1}}\otimes  \underline  {V_{2}}$ associated to the actions of $\widehat G$ on $V_{1}$ and $V_{2}$  commute with each other and the  diagonal action of $ \on{Gal}(\overline F/F)$ is the action  
 associated to the diagonal action of $\widehat G$ on
   $V_{1}\otimes V_{2}$.  
 This gives a structure as we want in a) because if $V$ is as above, $x\in V$, $\xi\in V^{*}$,  $f$ is the function on  $\widehat G$ defined as the matrix coefficient $f(g)=\langle \xi, g.x \rangle$, and $\gamma\in  \on{Gal}(\overline F/F)$ then we say that   $F_{f,\gamma}: \sigma\mapsto f(\sigma(\gamma))$, considered as a ``function on $\mathcal S$'',   acts on $H_{\{0\},\on{Reg}}$ by the composition   
 \begin{gather*} H_{\{0\},\on{Reg}} \xrightarrow{\on{Id}\otimes x} H_{\{0\},\on{Reg}} \otimes \underline V \xrightarrow{\gamma} H_{\{0\},\on{Reg}} \otimes \underline V 
 \xrightarrow{\on{Id}\otimes  \xi} H_{\{0\},\on{Reg}}. 
 \end{gather*}  
 Any function $f$ on $\widehat G$ can be written as such a matrix coefficient, 
 and the functions  $F_{f,\gamma}$ when $f$ and $\gamma$ vary are supposed to ``generate topologically all functions on $\mathcal S$''. The property above with $V_{1}$ and $V_{2}$ implies relations among the  $F_{f,\gamma}$, namely that 
 $F_{f,\gamma_{1}\gamma_{2}}=\sum_{\alpha }F_{f_{1}^{\alpha},\gamma_{1}}
 F_{f_{2}^{\alpha},\gamma_{2}}$ if the image of $f$ by comultiplication is
 $\sum_{\alpha} f_{1}^{\alpha}\otimes f_{2}^{\alpha}$. 
In \cite{zhu-notes} Xinwen Zhu 
gives an equivalent construction of the structure a). Structures a) and b) are compatible in the following sense: the conjugation $g F_{f,\gamma} g^{-1}$
of the action of $F_{f,\gamma}$ on $H_{\{0\},\on{Reg}} $ by the algebraic action of $g\in \widehat G$ is equal to the action of $F_{f^{g},\gamma} $ where 
$f^{g}(h)=f(g^{-1}hg)$. 
 
 The structures a) and b) give rise to a ``$\mathcal O$-module on the stack $\mathcal S/\widehat G$ of global Langlands parameters''  
 (such that the vector space of its ``global sections on $\mathcal S$'' is $H_{\{0\},\on{Reg}} $). 
For any morphism 
 $\sigma: \on{Gal}(\overline F/F)\to \widehat G(\overline{ \mathbb Q_{\ell}})$, 
  we want to define $\mathfrak A_{\sigma}$ as   the fiber of  this $\mathcal O$-module  at $\sigma$ (considered as a  ``$\overline{ \mathbb Q_{\ell}}$-point of $\mathcal S$ whose automorphism group in the stack $\mathcal S/\widehat G$ is $S_{\sigma}$'').  Rigorously we define $\mathfrak A_{\sigma}$ as the biggest quotient of $H_{\{0\},\on{Reg}}\otimes_{\mathbb Q_{\ell}}  \overline{ \mathbb Q_{\ell}}$ on which any function 
  $F_{f,\gamma}$ as above acts by multiplication by the scalar $f(\sigma(\gamma))$, and $S_{\sigma}$ acts on $\mathfrak A_{\sigma}$. If the heuristics  \eqref{heuri} is true it is the same as  $\mathfrak A_{\sigma}$ from the heuristics. 
  When $\sigma$ is elliptic (i.e. 
 when $S_{\sigma}/Z_{\widehat G}$ is finite), $\sigma$ is ``isolated in $\mathcal S$'' in the sense that it cannot be deformed (among 
 continuous morphisms whose composition with the abelianization of $\widehat G$ is of fixed  finite order) 
 and, as noticed by Xinwen Zhu,  heuristics  \eqref{heuri} is  true when we restrict on both sides to the parts lying over $\sigma$. In general due to deformation of some non elliptic $\sigma$ there could a priori be nilpotents, and for example we don't know how to prove 
   that $\mathcal B$ is reduced so we don't know how to prove the heuristics  \eqref{heuri}. 
      
We can see the  heuristics  \eqref{heuri}, and the structures a) and b) above,  as an illustration of the general idea that, 
 in a  spectral decomposition,  when  the points of the  spectrum  
naturally  have automorphism groups, the multiplicities should be  associated to representations of these groups.  By contrast the algebra $\mathcal B$ of excursion operators gives  the spectral decomposition with respect to the coarse quotient  associated to $\mathcal S/\widehat G$, where we forget the  automorphism groups $S_{\sigma}$. 
\end{rem}

 \begin{rem}
 The previous   remark makes sense although  $\mathcal S$ was not defined.  To define  a space like $\mathcal S$ rigorously  it may be necessary to consider continuous morphisms  $\sigma:\on{Gal}(\overline F/F)\to 
 \widehat G$ 
with coefficients in    $\mathbb Z_{\ell}$-algebras where $\ell$ is nilpotent (such $\sigma$ have  finite image), and $\mathcal S$ would be an  ind-scheme over $\on{Spf} \mathbb Z_{\ell}$. Then to define structure a) we would need      to consider 
$\on{Reg}$ with coefficients in $
\mathbb Z_{\ell}$, and,  for any representation $W$ of $\widehat G^{I}$ with coefficients   in $
\mathbb Z_{\ell}$,  to construct 
$H_{I,W}$ as a $
\mathbb Z_{\ell}$-module. 

 \end{rem}

\section{Local aspects: joint work with  Alain Genestier} 

  In   \cite{genestier-lafforgue}, Alain Genestier and I construct  
  the local parameterization up to  semisimplification   and the local-global compatibility.  
  
 Let $G$  be a  reductive group over a  local field  $K$ of equal characteristics. 
   We recall that the  Bernstein center of  $G(K)$ is defined, in two equivalent ways, as 
    \begin{itemize}
    \item  the center of the category of smooth representations of  $G(K)$, 
    \item  the algebra of central distributions on $G(K)$ acting as multipliers on the  algebra of locally constant  functions  with compact support. 
     \end{itemize}
  On  every  $\overline {{\mathbb Q}_{\ell}}$-linear   irreducible smooth representation of  $G(K)$, 
    the  Bernstein center acts by  a character.   
     
   The main result of  \cite{genestier-lafforgue}  associates to any   character $\nu$  of the  Bernstein center of $G(K)$ with values in $\overline{{\mathbb Q}_\ell}$  a local Langlands parameter $\sigma_{K}(\nu)$ up to   semisimplification , i.e.  (assuming  $G$ split to simplify) a conjugacy class of morphisms   $\on{Weil}(\overline K/K)\to \widehat G(\overline{{\mathbb Q}_\ell})$ defined over a finite extension of  ${\mathbb Q}_\ell$, continuous  and  semisimple.

We show in  \cite{genestier-lafforgue} the local-global compatibility up to  semisimplication, whose statement is the following. 
Let $X$  be a smooth  projective and geometrically irreducible curve over ${{\mathbb F}_q}$ and let $N$ be a level.  Then if  $\sigma$ is a  global Langlands parameter  and   if $\pi=\bigotimes \pi_{v}$  is an  irreducible representation of  $G(\mathbb A)$ such that  $\pi^{K_{N}}$ is non zero and appears in  $\mathfrak H_{\sigma}$ in the  decomposition \eqref{dec-nivau-N},   then for every  place $v$ de $X$ we have equality between 
  \begin{itemize}
  \item  the semisimplification of the  restriction  of  $\sigma$ to  $\on{Gal}(\overline{F_{v}}/F_{v})\subset \on{Gal}(\overline{F }/F )$,   
 \item  the semisimple local parameter 
 $\sigma_{K}(\nu)$ where $\nu$ is the character of the Bernstein center by which it acts on   the   irreducible  smooth representation 
  $ \pi_{v} $ of  $G(K)$.    
  \end{itemize}

  We use nearby cycles on arbitrary bases (Deligne, Laumon, Gabber, Illusie, Orgogozo),  which are defined on  oriented products of toposes and whose properties are proven in \cite{orgogozo}  (see also \cite{illusie} for an excellent  survey). 
    Technically we show that if all the 
    $\gamma_{i}$ are in   $\on{Gal}(\overline{F_{v}}/F_{v})\subset \on{Gal}(\overline F/F)$
    then the global excursion operator   $    S_{I,f,(\gamma_{i})_{i\in I}} \in \on{End}\big( C_{c}^{\mathrm{cusp}}(\on{Bun}_{G,N}({{\mathbb F}_q})/\Xi,{\mathbb Q}_\ell)\big)$ acts by  multiplication by an element 
 $\mathfrak z_{I,f,(\gamma_{i})_{i\in I}}$ of the $\ell$-adic completion of the  Bernstein center of  $G(F_{v})$ which depends only on the local data at  $v$. We construct 
  $\mathfrak z_{I,f,(\gamma_{i})_{i\in I}}$ using stacks of ``restricted shtukas'', which are analogues of truncated Barsotti-Tate groups. 
      
\begin{rem}
In the case of  $GL_{r}$ the local  correspondence was   known by Laumon-Rapoport-Stuhler \cite{laumon-rapoport-stuhler} and the   local-global compatibility  (without semisimplification) was proven in \cite{laurent-inventiones}. 
Badulescu and Henniart explained us that in general 
we cannot hope more that the   local-global compatibility up to   semisimplication. 
\end{rem}

\section{Independence on $\ell$}

 Grothendieck motives (over a given field) form a  $\overline {\mathbb Q}$-linear category and unify the $\ell$-adic cohomologies  (of varieties over this field) for different  $\ell$: a motive is  ``a factor in  a universal cohomology of a variety''.   We consider here motives over $F$. 
We conjecture that the decomposition $$C_{c}^{\mathrm{cusp}}(\on{Bun}_{G}({{\mathbb F}_q})/\Xi,\overline{{\mathbb Q}_\ell})=\bigoplus_{\sigma} \mathfrak H_{\sigma}$$ we have constructed  is defined over  $\overline{\mathbb Q}$ (instead of $\overline{{\mathbb Q}_\ell}$), 
  indexed by  motivic     Langlands parameters  $\sigma$, 
   and independent on  $\ell$. This conjecture seems out of reach for the moment.

  The  notion of 
 motivic Langlands parameter is clear if we admit the standard conjectures. 
 A  motivic Langlands parameter defined over $\overline{\mathbb Q}$
 would give rise to a ``compatible'' family of morphisms $\sigma_{\ell, \iota}:\on{Gal}(\overline F/F)\to \widehat G(\overline{{\mathbb Q}_\ell})$
 for any $\ell$ not dividing $q$ and any embedding $\iota:\overline{\mathbb Q}\hookrightarrow  \overline{{\mathbb Q}_\ell}$. 
 When $G=GL_{r}$, the condition of compatibility is straightforward 
 (the traces of the Frobenius elements should belong to $\overline{\mathbb Q}$ and be the same for all $\ell$ and $\iota$) 
 and the fact that any irreducible representation (with determinant of finite order) for some $\ell$ belongs to such a family
 (and has therefore ``compagnons'' for other $\ell$ and $\iota$) 
 was proven as a consequence of the Langlands correspondence in \cite{laurent-inventiones}. It was generalized in the two following independent directions  
 \begin{itemize}
 \item Abe \cite{abe} used the crystalline cohomology of stacks of shtukas
 to construct  crystalline compagnons,   
 \item  when $F$ is replaced by the field of rational functions of  any smooth variety over ${{\mathbb F}_q}$, Deligne proved that the traces of Frobenius elements belong to a finite extension of ${\mathbb Q}$ and Drinfeld 
 constructed these compatible families  \cite{deligne-finitude, drinfeld-del-conj, esnault-deligne}.  
 \end{itemize}
 
 For a general reductive group $G$ the notion of compatible family is subtle 
 (because the obvious condition on the conjugacy classes  of the Frobenius elements is 
 not sufficient).  In \cite{drinfeld-pro-completion} Drinfeld 
 gave the right conditions to define compatible families and proved that any continuous semisimple morphism 
 $\on{Gal}(\overline F/F)\to \widehat G(\overline{{\mathbb Q}_\ell})$  factorizing through 
            $\pi_{1}(U, \overline \eta)$  for some open dense $U\subset X$ (and such that the Zariski closure of its image is semisimple) belongs to  a unique compatible family.

   \section{Conjectures on Arthur parameters}
   
   We hope that all 
  global Langlands parameters  $\sigma$ which appear in this  decomposition come from elliptic  Arthur parameters, i.e. conjugacy  classes of continuous semisimple morphisms $\on{Gal}(\overline F/F)\times SL_{2}(\overline{{\mathbb Q}_\ell})\to \widehat G(\overline{{\mathbb Q}_\ell})$ whose  centralizer is finite  modulo the center of  $\widehat G$. This  $SL_{2}$ should be related to the  Lefschetz  $SL_{2}$ acting on the  intersection cohomology of  compactifications of stacks of shtukas. 
We even hope a  parameterization of the vector space of discrete  automorphic forms (and not only cuspidal ones)  indexed by elliptic Arthur parameters.  

Moreover we expect that generic cuspidal automorphic forms appear exactly in $\mathfrak H_{\sigma}$ such that $\sigma$ is elliptic as a Langlands parameter
(i.e.  that it comes from an elliptic Arthur parameter with trivial $SL_{2}$ action). 
This would imply the Ramanujan-Petersson conjecture (an archimedean estimate on Hecke eigenvalues).

  By \cite{drinfeld-kedlaya} the conjectures above would also imply $p$-adic estimates on Hecke eigenvalues which would sharper than those in \cite{estimees-SMF}.

  \section{Recent works on the Langlands program for function fields in relation with shtukas}\label{other-works}

In    \cite{boeckle-harris...} G.  B\"ockle, M. Harris, C. Khare, and J. Thorne
   apply the results explained in this text together with Taylor-Wiles methods to prove (in the split and everywhere unramified situation) the potential automorphy of all Langlands parameters with Zariski-dense image. Thus they prove a weak form of the ``Galois to automorphic'' direction.    
   
In   
   \cite{Yun-Zhang}     Zhiwei Yun and Wei Zhang proved  analogues of the Gross-Zagier formula, namely equality between the intersection numbers of some algebraic cycles in stacks of     shtukas and  special values of derivatives of L-functions (of arbitrary order equal  to the number of legs). 
   
    In    \cite{zhu} 
Liang Xiao and  Xinwen Zhu construct algebraic cycles in special fibers of Shimura varieties. Their construction was inspired by the case of the stacks of shtukas and is  already new in this case (it gives  a conceptual setting for the Eichler-Shimura relations   used in \cite{coh}). 

\section{Geometric Langlands program}

The results   explained above are based on the geometric Satake equivalence  \cite{mv}, and are inspired by the factorization structures studied by   Beilinson and Drinfeld
 \cite{chiral}.  
The geometric Langlands program was pioneered by Drinfeld \cite{drinfeld-geom} and Laumon  \cite{laumon-duke}, and then developped itself in  two variants, which we will discuss in turn. 

\subsection{Geometric Langlands program for $\ell$-adic sheaves}
 Let $X$ be a smooth projective curve over an algebraically closed field of characteristic different from $\ell$. 
 
 For any representation $W$ of $\widehat G^{I}$  the  Hecke functor 
     $$\phi_{I,W}:D^{b}_{c}(\on{Bun}_{G},\overline{{\mathbb Q}_\ell})\to  D^{b}_{c}(\on{Bun}_{G}\times X^{I}, \overline{{\mathbb Q}_\ell})$$      
    is given by  
     $$\phi_{I,W}(\mathcal F)=q_{1,!}\big(q_{0}^{*}(\mathcal F)\otimes \mathcal F_{I,W}\big)$$
 where $\on{Bun}_{G}\xleftarrow{q_{0}}\on{Hecke}_{I}\xrightarrow{q_{1}}\on{Bun}_{G}\times X^{I}$ is the Hecke   correspondence classifying modifications of a $G$-bundle at the $x_{i}$, and  $\mathcal F_{I,W}$ is defined as   the inverse image of  
 $\mathcal S_{I,W}$ by the natural formally smooth morphism   
 $\on{Hecke}_{I}\to \mathcal M_{I}$.  
 
  Let $\sigma$  be a  $\widehat G$-local system on  $X$. Then  $\mathcal F\in D^{b}_{c}(\on{Bun}_{G},\overline{{\mathbb Q}_\ell})$ is said to be an eigensheaf for   $\sigma$ if we have, for any finite set $I$ and any representation $W$ of 
  $(\widehat G)^{I}$,  an isomorphism 
  $\phi_{I,W}(\mathcal F)\overset {\thicksim}{\to} \mathcal F\boxtimes W_{\sigma}$,   functorial in   $W$ and compatible to exterior tensor products and fusion. The    conjecture of the geometric  Langlands program claims the existence of an  $\sigma$-eigensheaf  $\mathcal F$  (it should also satisfy a  Whittaker normalization condition which in particular prevents it to be $0$). For $G=GL_{r}$ this conjecture was proven by Frenkel, Gaitsgory,  Vilonen in \cite{fgv, ga-vanishing} 

When $X$, $\on{Bun}_{G}$, $\sigma$ and $\mathcal F$ are defined over ${{\mathbb F}_q}$ (instead of $\overline  {{\mathbb F}_q}$), a construction of Braverman and Varshavsky \cite{brav-var} produces subspaces of cohomology classes in the stacks of shtukas and this allows to show that the 
function given by the trace of Frobenius on $\mathcal F$ belongs to 
 the factor $\mathfrak H_{\sigma}$ of  decomposition \eqref{dec-intro-thm-ppal}, as explained in section 15 of \cite{coh}.

 The $\ell$-adic setting is truly a geometrization of automorphic forms over function fields, and many constructions were geometrized: 
Braverman and Gaitsgory geometrized Eisenstein series \cite{bg-eis}, and Lysenko geometrized in particular  Rankin-Selberg integrals \cite{sergey-RS}, theta correspondences
 \cite{sergey-theta,sergey-theta-SO-Sp,even-orth}, and several constructions for  metaplectic groups \cite{lysenko-twisted-torus,lysenko-whitt-meta}.

 \subsection{Geometric Langlands program for $D$-modules}
 Now let  $X$ be a smooth projective curve over an algebraically closed field of characteristic $0$. 
A feature of the setting of $D$-modules is that one can upgrade the statement of Langlands correspondence 
to a conjecture about an equivalence between categories on the geometric and spectral sides, respectively.
See \cite{dennis-laumon} for a precise statement of the conjecture  and \cite{dennis-bki} for a survey of recent progress. Such an equivalence can in principle make sense 
due to the fact that Galois representations into $\widehat G$, instead of being taken individually, now
form an algebraic stack $\on{LocSys}_{\widehat{G}}$ classifying  $\widehat{G}$-local systems, i.e. $\widehat{G}$-bundles with connection (by contrast one does not have such an algebraic stack in the $\ell$-adic setting).

On the geometric side, one considers the 
derived 
category of D-modules on $\on{Bun}_G$, or rather a stable $\infty$-category enhancing it. It is denoted   $\on{D-mod}(\on{Bun}_G)$  and is defined and studied in \cite{dennis-drinfeld-Dmod}. 
The category on the spectral side is a certain modification of $\on{QCoh}(\on{LocSys}_{\widehat{G}})$,
the (derived or rather $\infty$-) category of quasi-coherent sheaves on the stack $\on{LocSys}_{\widehat{G}}$. The modification in question
is denoted $\on{IndCoh}_{\on{Nilp}}(\on{LocSys}_{\widehat{G}})$, and it has to do with the fact that $\on{LocSys}_{\widehat{G}}$ is not smooth, but
rather \emph{quasi-smooth} (a.k.a. derived locally complete intersection). The difference between $\on{IndCoh}_{\on{Nilp}}(\on{LocSys}_{\widehat{G}})$
and $\on{QCoh}(\on{LocSys}_{\widehat{G}})$ is measured by \emph{singular support} of coherent sheaves, a theory developed in \cite{AG}. The introduction of $\on{Nilp}$ in \cite{AG} was motivated by the case of $\mathbb P^1$ \cite{cas-P1} and the study of the singular support  of the geometric Eisenstein series.  
In terms of Langlands correspondence, this singular support may also be seen as accounting for \emph{Arthur parameters}. More precisely 
the singularities of $\on{LocSys}_{\widehat{G}}$ are controlled by a 
stack 
$\on{Sing}(\on{LocSys}_{\widehat{G}})$ over $ \on{LocSys}_{\widehat{G}}$ whose fiber over a point $\sigma$ is the $H^{-1}$ of the cotangent complex at $\sigma$, equal to $H^{2}_{dR}(X, \widehat {\mathfrak g}_{\sigma})^{*}\simeq H^{0}_{dR}(X, \widehat {\mathfrak g}_{\sigma}^{*})
\simeq  H^{0}_{dR}(X, \widehat {\mathfrak g}_{\sigma}) $
where the first isomorphism is Poincaré duality and the second depends on the 
 choice of a non-degenerate $ad$-invariant symmetric bilinear form on $\widehat {\mathfrak g}$.  Therefore $\on{Sing}(\on{LocSys}_{\widehat{G}})$ is identified to  the stack classifying $(\sigma, A)$, with 
$\sigma\in \on{LocSys}_{\widehat{G}}$ and $A$ an horizontal section of the local system $ \widehat {\mathfrak g}_{\sigma}$
associated to $\sigma$ with the adjoint representation of $\widehat G$.  
Then $\on{Nilp}$ is   the cone of $\on{Sing}(\on{LocSys}_{\widehat{G}})$  defined by the condition that $A$ is nilpotent. 
By the Jacobson-Morozov theorem, any such $A$ is the nilpotent element associated to a morphism of $SL_{2}$ to the centralizer of $\sigma$ in $\widehat G$, i.e. it comes from an Arthur parameter. 
The singular support  
of a coherent sheaf on $\on{LocSys}_{\widehat{G}}$ is a closed substack in 
$\on{Sing}(\on{LocSys}_{\widehat{G}})$. 
 The category $\on{IndCoh}_{\on{Nilp}}(\on{LocSys}_{\widehat{G}})$
 (compared to $\on{QCoh}(\on{LocSys}_{\widehat{G}})$) corresponds to  the condition that the singular support of coherent sheaves has to be included in $\on{Nilp}$ (compared to the zero section where $A=0$). 
The main conjecture is that there is an equivalence of categories 
\begin{gather}\label{eq-cat-conj}\on{D-mod}(\on{Bun}_G)\simeq \on{IndCoh}_{\on{Nilp}}(\on{LocSys}_{\widehat{G}}).\end{gather}

Something weaker is known: 
by \cite{dennis-laumon},  $\on{D-mod}(\on{Bun}_G)$ ``lives" over $\on{LocSys}_{\widehat{G}}$ in the sense  that
$\on{QCoh}(\on{LocSys}_{\widehat{G}})$, viewed as a monoidal category, acts naturally on $\on{D-mod}(\on{Bun}_G)$.
Note that $\on{QCoh}(\on{LocSys}_{\widehat{G}})$ acts on $\on{IndCoh}_{\on{Nilp}}(\on{LocSys}_{\widehat{G}})$
(one can tensor a coherent complex by a perfect one and obtain a new coherent complex) and the conjectured equivalence \eqref{eq-cat-conj} should be compatible with the actions of $\on{QCoh}(\on{LocSys}_{\widehat{G}})$ on both sides.

Theorem~\ref{intro-thm-ppal} (refined in remark~\ref{rem-drinfeld-construction})  can be considered as an arithmetic analogue of  the fact that 
$\on{D-mod}(\on{Bun}_G)$ ``lives" over $\on{LocSys}_{\widehat{G}}$  (curiously, due to the lack of an $\ell$-adic analogue of  $\on{LocSys}_{\widehat{G}}$, that result does not have an analogue  
 in the $\ell$-adic geometric Langlands program, even if the vanishing conjecture proven by  Gaitsgory \cite{ga-vanishing} goes in this direction). And the fact that Arthur multiplicities formula is still unproven in general is parallel to the fact that the equivalence \eqref{eq-cat-conj}  is still unproven.

When $G=T$ is a torus, there is no difference between $\on{QCoh}(\on{LocSys}_{\widehat T})$ and $\on{IndCoh}_{\on{Nilp}}(\on{LocSys}_{\widehat T})$. 
In this case, the desired equivalence $\on{QCoh}(\on{LocSys}_{\widehat T})\simeq \on{D-mod}(\on{Bun}_T)$ is a theorem, due to  Laumon \cite{fourier-gen}. 

The formulation of the geometric Langlands correspondence as an equivalence of categories 
\eqref{eq-cat-conj}, and even more the proofs of the results,  rely on substantial developments
in the technology, most of which had to do with the incorporation of the tools of higher category
theory and higher algebra, developed by J.~Lurie in \cite{Lu1,Lu2}.  Some of the key constructions
use categories of D-modules and quasi-coherent sheaves on algebro-geometric objects more
general than algebraic stacks (a typical example is the moduli space of $G$-bundles on $X$ equipped
with a reduction to a subgroup \emph{at the generic point of $X$}).

\subsection{Work of Gaitsgory and Lurie on Weil's conjecture on Tamagawa numbers over function fields}

In \cite{GL1,GL2} (see also \cite{Atiyah-Bott})  Gaitsgory and Lurie compute the cohomology with coefficients in $ \mathbb Z_{\ell}$ of the stack $\on{Bun}_{G}$ when $X$ is any smooth projective curve over an algebraically closed  field of characteristic other than $\ell$, and $G$ is a smooth  affine group scheme over X with connected
fibers, whose generic fiber is semisimple and simply connected. 
They use in particular a remarkable geometric ingredient, belonging to the same framework 
of factorization structures \cite{chiral} (which comes from conformal field theory) 
as the geometric Satake equivalence. 
The Ran space of $X$ is, loosely speaking,  the  prestack classifying non-empty finite subsets $Z$ of $X$. 
The affine grassmannian  $\on{Gr}_{Ran}$  is the  prestack over the Ran space classifying such a $Z$, 
a $G$-bundle $\mathcal G$ on $X$, and  a trivialization $\alpha$ of $\mathcal G$ on $X\setminus Z$. Then the remarkable geometric ingredient is that the obvious morphism $\on{Gr}_{Ran}\to \on{Bun}_{G}, 
(Z, \mathcal G, \alpha)\mapsto \mathcal G$ has contractible fibers in some sense and gives an isomorphism on homology.  Note that when $k=\mathbb C$ and $G$ is constant on the curve, their formula implies the well-known Atiyah-Bott formula for the cohomology of $\on{Bun}_{G}$, whose usual proof is by analytic means. 

Now assume that the curve $X$ is over $\mathbb F_{q}$. By the Grothendieck-Lefschetz trace formula their computation of the cohomology of $\on{Bun}_{G}$ over $\overline{\mathbb F_{q}}$ gives a formula for  $|\on{Bun}_G(\mathbb F_q)|$, 
the number of $\mathbb  F_q$-points on the stack $\on{Bun}_G$. Note that
since $\on{Bun}_G$ is a stack, each isomorphism class $y$ of points is weighted by $\frac{1}{\on{Aut}_y(\mathbb F_q)}$, where $\on{Aut}_y$ is the
algebraic group of automorphisms of $y$, and $\on{Aut}_y(\mathbb  F_q)$ is the finite group of its $\mathbb F_q$-points. 
Although the set of isomorphism classes $y$ of points is infinite,   the weighted sum   converges. Gaitsgory and Lurie easily  reinterpret $|\on{Bun}_G(\mathbb F_q)|$ as the volume 
(with respect to some measure)   
of the quotient $G(\mathbb A)/G(F)$ (where $F$ is the function field of $X$ and 
$\mathbb A$ is its ring of adèles) and prove in this way, in the case of function fields,  a formula for the volume of $G(\mathbb A)/G(F)$ as a product of local factors at all places. This formula, called the  Tamagawa number formula,  had been conjectured by  Weil for any global field $F$. 

Over number fields $\on{Bun}_G$ does not make sense, only the conjecture of Weil on the Tamagawa number formula remains and it had been proven by Kottwitz after earlier works of Langlands and Lai by completely different methods (residues of Eisenstein series and trace formulas).

 \section{Homage to Alexandre Grothendieck (1928-2014)}

Modern algebraic geometry 
was built  by
  Grothendieck, together with his students, 
in the realm of categories: 
 functorial definition of schemes and stacks, 
Quot construction  for  $\on{Bun}_{G}$, tannakian formalism, topos,  étale cohomology, 
motives. 
 His vision 
 of topos and motives  already had tremendous consequences and 
 others are certainly yet to come. 
He also had a strong influence outside of his school,  
 as testified by the rise of higher categories and the work of  
Beilinson, Drinfeld, Gaitsgory, Kontsevich, Lurie, Voevodsky (who, sadly, passed away recently) and many others. 
 He changed not only mathematics, but also the way we think about it.

\end{document}